\documentclass[11pt]{article}

\usepackage{amsmath,amsfonts,amssymb,amsthm}

\usepackage[all]{xy}

\newtheorem{thm}{Theorem}[section]

\newtheorem{dfn}[thm]{Definition}

\newtheorem{lem}[thm]{Lemma}

\newtheorem{prop}[thm]{Proposition}

\newtheorem{remark}[thm]{Remarks}

\newtheorem{cor}[thm]{Corollary}

\newtheorem{ex}[thm]{Example}

\newtheorem{question}[thm]{Question}

\def\sq{{\scriptscriptstyle \square}}

\def\bq{\begin{question}}
\def\bt{\begin{thm}}
\def\bp{\begin{prop}}
\def\blem{\begin{lem}}
\def\bd{\begin{dfn}}
\def\br{\begin{remark}}
\def\bc{\begin{cor}}
\def\bex{\begin{ex}}
\def\beqs{\begin{eqnarray*}}
\def\beq{\begin{eqnarray}}
\def\bi{\begin{itemize}}

\def\eq{\end{question}}
\def\et{\end{thm}}
\def\ep{\end{prop}}
\def\elem{\end{lem}}
\def\ed{\end{dfn}}
\def\er{\end{remark}}
\def\ec{\end{cor}}
\def\eex{\end{ex}}
\def\eeqs{\end{eqnarray*}}
\def\eeq{\end{eqnarray}}
\def\ei{\end{itemize}}

\def\wot{\widehat{\otimes}}

\def\ds{\displaystyle}

\def\c{\cdot}

\def\ov{\overline}
\def\r{\rangle}
\def\l{\langle}
\def\prep{\{ \pi, {\cal H} \} }

\def\H{{\cal H}}

\def\B{B({\cal H})}

\def\T{\mathbb{T}}

\def\w*{{\rm wk}^*-{\rm wk}^*}

\def\ra{\rightarrow}

\def\C{\mathbb{C}}

\def\xpe{\xi*_\pi \eta}
\def\epx{\eta*_\pi \xi}

\def\F{{\cal F}}

\def\T{{\cal V}}

\def\wot{\widehat{\otimes}}

\def\TA{{\cal T}(A^*)}
\def\T{{\cal T}}

\def\SA{{\cal S}(A^*)}
\def\S{{\cal S}}
\def\E{{\cal E}}

\def\FqA{{\cal F}_q(A^*)}
\def\FQA{{\cal F}_Q(A^*)}
\def\F{{\cal F}}

\def\EqA{{\cal E}_q(A^*)}
\def\EQA{{\cal E}_Q(A^*)}

\def\A{{\mathfrak A}}
\def\B{{\mathfrak B}}

\def\etaa{\eta_{\A}}
\def\etab{\eta_{\B}}

\def\oH{\ov{\H}}

\def\bg{\breve{\Gamma}} 

\begin{document}
\title{Fourier spaces  and  Completely Isometric Representations of Arens Product Algebras}
\author{Ross  Stokke \footnote{
This research was partially supported   by an NSERC grant. }    }
\date{}
\maketitle

\vspace{-.1in}

\begin{abstract}{\small   Motivated by the definition of a semigroup compactification of a locally compact group and a large collection of examples, we introduce the notion of an (operator) ``homogeneous left dual Banach algebra" (HLDBA) over a (completely contractive) Banach algebra  $A$.  We prove a Gelfand-type representation theorem showing that every HLDBA over $A$ has a concrete realization as an (operator) homogeneous left Arens product algebra:  the dual of a subspace of $A^*$ with a compatible (matrix) norm and a type of left Arens product $\sq$.  Examples include all left Arens product algebras over $A$, but also -- when $A$ is the group algebra of a locally compact group -- the dual of its Fourier algebra.   Beginning with any (completely) contractive (operator) $A$-module action $Q$ on a space $X$, we introduce the (operator) Fourier space $(\FQA, \| \c \|_Q)$ and 
 prove that $(\FQA^*, \sq)$ is the unique (operator) HLDBA  over $A$ for which there is a weak$^*$-continuous completely isometric representation as completely bounded operators on $X^*$  extending the dual module representation. 
   Applying our theory to several examples of (completely contractive) Banach algebras $A$ and module operations, we provide new characterizations of 
  familiar HLDBAs over $A$ and we recover -- and often extend -- some  
  (completely) isometric representation theorems  concerning these HLDBAs.  

\smallskip

\noindent{\em Primary MSC codes:}  47L10, 47L25,  43A20,  43A30,   46H15,  46H25  
   
 }
\end{abstract}

\section{Introduction}

Many of the most well-studied and basic objects associated with a locally compact group $G$ -- more generally a locally compact quantum group -- are introverted subspaces  of $L^1(G)^*=L^\infty(G)$ and their dual spaces under an Arens product: examples include the introverted space of continuous functions vanishing at infinity, $C_0(G)$ (its dual with Arens product is the measure algebra $M(G)$ with convolution product); the introverted space of continuous almost periodic functions on $G$, $AP(G)$;  the introverted space of continuous Eberlein functions on $G$, $E(G)$; the introverted space of continuous weakly almost periodic functions on $G$, $WAP(G)$; the left introverted space of left uniformly continuous  functions on $G$, $LUC(G)$; and $L^\infty(G)$. A small sample of papers in which the duals of these and other  spaces are studied as  left Arens product algebras  is \cite{Are, Dal-Lau, Fil-Mon-Neu, Hu-Neu-Rua, Lau1, Lau2, Neu2, Neu-Rua-Spr, Spr-Sto}. 

In general, if $\SA$ is a closed left introverted subspace of the dual space, $A^*$, of a Banach algebra $A$, then --  with its left Arens product $\sq$ --  $\A= \SA^*$ is a Banach algebra such that 
\bi \item[(i)] multiplication is separately weak$^*$-continuous with respect to every fixed right variable and weak$^*$-dense ``topological centre" $Z_t(\A)$ ;  and 
\item[(ii)] there is a continuous homomorphism  $\eta_\A: A \ra \A$   mapping into $Z_t(\A)$  such that 
\item[(iii)] the image of the unit ball of $A$ under $\eta_\A$ is weak$^*$-dense in the unit ball of $\A$. 
\ei
 We call $\A = \SA^*$ a left Arens product algebra over $A$. In \cite{Sto1} we called a pair $(\A, \eta_\A)$ a left dual Banach algebra (LDBA) over $A$ when $\A$ is a Banach algebra and a dual space such that   properties (i), (ii) and (iii) are satisfied. Motivated by the fundamental theorem from semigroup compactification theory stating that  every right topological semigroup compactification of a locally compact group is a Gelfand compactification, we proved that, up to equivalence, every LDBA over $A$ is a left Arens product algebra over $A$ \cite[Theorem 3.3] {Sto1}. 
 
 In addition to the class of examples provided by the left Arens product algebras over  $A$, the reader will likely be aware of several examples of pairs $(\A, \eta_\A)$ where $\A$ is a Banach algebra and a dual Banach space such that properties (i), (ii) and  -- a weaker version of statement (iii) -- 
 \bi \item[(iii')] the image of $A$ under $\eta_\A$ is weak$^*$-dense in $\A$  
 \ei 
 are satisfied.  For instance, if $A$ is any involutive Banach algebra and $\eta_\A: A \ra B(\H)$ is a $*$-representation of $A$ on a Hilbert space $\H$, letting $\A$ denote the von Neumann subalgebra of $B(\H)$ generated by $\eta_\A(A)$, the pair  $(\A, \eta_\A)$ satisfies properties (i),(ii) and (iii'). We will call any such pair, ($\A, \eta_\A)$,  a homogeneous left dual Banach algebra (HLDBA) over $A$. 
 
 Our Gelfand-type representation theorem from \cite{Sto1} suggests that any HLDBA over $A$ looks something like a left Arens product algebra, and in Section 3 we introduce the notion of a left introverted homogeneous subspace of $A^*$, $(\S(A^*), \| \c \|_\S)$ and an Arens-type product $\sq$ on $\SA^*$ such that the pair $(\SA^*, \eta_\S)$ is a HLDBA over $A$, where $\eta_\S$ is defined by  $\l \eta_\S(a), \phi\r = \l \phi, a\r$. We call $\SA^*$ a homogeneous left Arens product algebra over $A$ and we prove a new  Gelfand-type representation theorem stating that every HLDBA over $A$ is equivalent to a homogeneous left Arens product algebra over $A$. By introducing a notion of subdirect product for HLDBAs over $A$,  we also show that with respect to a natural ordering $\leq$, the partially ordered set $({\cal HLD}(A), \leq)$ of all HLDBAs over $A$ is a complete lattice. The  notion of subdirect product in the category of semigroup compactifications  is known to provide an efficient method of constructing universal semigroup compactifications \cite{Ber-Jun-Mil}. The results in Section 3 are proved in the setting of operator spaces. 
 
 A particular example of a left introverted homogeneous subspace of $L^1(G)^*$ is $(A(G), \| \c \|_A)$, the Fourier algebra of $G$ with its Fourier norm. Thus, $(A(G)^*, \eta_A)$ with the Arens-type product $\sq$ defined in Section 3 is an example of a HLDBA over $L^1(G)$. In this case, we can identify  $A(G)^*$ with  $VN(G)$, the von Neumann subalgebra of $B(L^2(G))$ generated by the left regular representation $\{\lambda_2, L^2(G)\}$ of $L^1(G)$.  In fact,  the identifying map $\Phi: A(G)^* \ra VN(G)$ is a weak$^*$-continuous isometric algebra isomorphism with respect to $\sq$ such that $ \Phi \circ \eta_A = \lambda_2$. In the language of this paper, $(A(G)^*, \eta_A)$ and $(VN(G), \lambda_2)$ are equivalent HLDBAs over $L^1(G)$, and $\Phi$ is a weak$^*$-continuous isometric representation of the Banach algebra $(A(G)^*, \sq)$ on $L^2(G)$.  The pair $(A(G)^*, \eta_A)$ is one of many known examples of  an HLDBA over $A$  that  can be isometrically represented as an algebra of operators on a Banach space $E$. 
 
 The main constructions and results of this paper are found in Section 4. Beginning with a completely contractive Banach algebra $A$ and a completely contractive right operator $A$-module action $Q:X\times A \ra X$, in Section 4.1 we construct the operator Fourier space $\FQA$ together with an associated operator space matrix norm $\| \c \|_Q$ such that $(\FQA, \| \c \|_Q)$ is an operator left introverted homogeneous subspace of $A^*$ and   $(\FQA^*, \sq)$ has a weak$^*$-continuous completely isometric representation as completely bounded mappings on $X^*$.  We characterize $(\FQA^*, \eta_Q)$ as the unique HLDBA over $A$ with this property. Within the complete lattice $({\cal HLD}(A), \leq)$, we characterize $(\FQA^*, \eta_Q)$ as the minimum HLDBA over $A$ with a weak$^*$-continuous completely contractive representation in $CB(X^*)$  that extends the dual $A$-module action on $X^*$ determined by $Q$.    In Section 4.2, we introduce the corresponding construction of the Fourier space $\FqA$ and Fourier norm $\| \c \|_q$ in the category of Banach spaces. 
 
 In Section 5, the theory developed in Sections 3 and 4 is applied to different module operations $X \times A \ra X$. We provide new characterizations of some well-studied (operator) homogeneous left Arens product algebras over $A$ and we 
  recover -- and often extend -- some  familiar results concerning weak$^*$-continuous (completely) isometric representations of these Banach algebras:  
 \bi \item When $A$ is an involutive Banach algebra and $\{\pi, \H\}$ is a $*$-representation of $A$ on a Hilbert space $\H$, in Section 5.1 we introduce the Fourier space $(\F_\pi(A^*), \| \c \|_\pi)$ and show that $(\F_\pi(A^*), \sq)$ is a $W^*$-algebra that can be identified with $VN_\pi$, the von Neumann subalgebra of $B(\H)$ generated by $\pi$, via a weak$^*$-continuous isometric $*$-isomorphism that extends $\pi$. When $A = L^1(G)$ and $\prep$ is a continuous unitary representation of $G$,  $(\F_\pi(A^*), \| \c \|_\pi)$ is the Arsac-Eymard Fourier space $(A_\pi, \| \c \|_{B(G)})$, where $\| \c \|_{B(G)}$ is the Fourier-Stieltjes algebra norm. Taking $\prep = \{ \lambda_2, L^2(G)\}$, $\F_\pi(A^*)$ is thus the Fourier algebra $A(G)$ and we recover the identification of $A(G)^*$ with $VN(G)$ and, further, we are able to identify the product on $VN(G)$ with $\sq$. 
 \item In Section 5.2, we recover the  Fig\`{a}-Talamanca--Herz spaces $A_p(G)$  and the identification of $A_p(G)^*$ with the operator subalgebra $PM_{p'}(G)$ of $B(L^{p'}(G))$ -- the algebra of $p'$-pseudomeasures. 
 \item Given a completely contractive Banach algebra $A$ with a contractive right approximate identity, in Section 5.3 we recognize $LUC(A^*)$, the left introverted subspace of $A^*$ comprised of left uniformly continuous functionals on $A$, as the Fourier space $\FQA$ associated with the right module action $Q$ of $A$ on itself. This allows us to  identify the completely contractive Banach algebra $(LUC(A^*)^*, \sq)$  with the weak$^*$-closed subalgebra $CB_A(A^*)$ of completely bounded right $A$-module maps on $A^*$ via a weak$^*$-homeomorphic completely isometric algebra isomorphism. This  includes results that were pioneered by Curtis and Fig\`{a}-Talamanca in \cite{Cur-Fig},   Lau in \cite{Lau1, Lau2}, and extended in \cite{Hu-Neu-Rua, Neu-Rua-Spr}, for various examples of $A$. 
 \item   By recognizing $LUC(G)^*$ as the Fourier space $\F_Q(L^1(G)^*)$  associated to a natural $L^1(G)$-module action on the trace class operators on $L^2(G)$, in Section 5.4 we recover a result of M. Neufang that provides a completely isometric representation of $LUC(G)^*$ as completely bounded mappings on $B(L^2(G))$ \cite{Neu1, Neu2}. Moreover, we show that $LUC(G)^*$ is characterized within the set (up to equivalence) of HLDBAs over $A$ by the existence of such a representation.   
 \ei   
 
\section{Preliminary definitions and results}

Let $A$ be a (completely contractive) Banach algebra. If $X$ is a Banach space, we will often use the  notation $x$,  $x^*$ and $x^{**}$  for elements in $X$, the dual $X^*$ and the bidual $X^{**}$ without explanation; an element $a$ or $b$ will always belong to the Banach algebra $A$. 
By a  (right, left, or bi-) $A$-module we will mean a Banach $A$-module.  A right $A$-module $X$ is \it essential \rm if the  closed linear span of $ X \c A $ is all of $X$ and $X$ is \it neo-unital \rm when $X =  X \c A$. By the Cohen factorization theorem \cite[Theorem 32.22]{HRII}, when $A$ has a  right bounded approximate identity, the two concepts are the same.    When $X$ is a right (left)  $A$-module, its dual space $X^*$ becomes a \it dual left (right) $A$-module \rm with products given by 
$$\l a\c x^*, x \r = \l x^*, x \c a \r  \quad \quad ( \l x^* \c a , x \r = \l x^*, a \c x \r).$$ In particular, $A^*$ will always be viewed as  dual Banach $A$-module with respect to the operations
$$ \l a \c \phi, b \r = \l \phi, ba\r  \ \ {\rm and} \ \   \l \phi \c a , b \r = \l \phi , a b  \r \quad \quad  \phi \in A^*.$$
Unless the codomain is a scalar field, all maps between normed linear spaces  are assumed to be linear and (norm)  continuous. 

We denote the operator space projective tensor product \cite[Ch. 7]{Eff-Rua} of operator spaces $X$ and $Y$ by $X \wot Y$ and will use the notation  $X \otimes^\gamma Y$ to denote a Banach space projective tensor product.  To aid the reader, we will typically use upper-case script -- $P, \ Q$, etc. --   to denote maps employed when working in the category of operator spaces and lower-case script  -- $p, \ q$, etc. --  when working in the category of Banach spaces.

Let $m: X \times Y \ra Z: (x,y) \mapsto x \c y $ be a bounded bilinear map. The (first)  Arens transpose of $m$  is the bounded bilinear map $$m': Z^* \times X \ra Y^*: (z^*, x) \mapsto z^* \c x$$ where $$\l m'(z^*,x), y \r = \l z^*, m(x,y) \r; \ \ {\rm or} \ \ \l z^* \c x, y \r = \l z^*, x \c y \r $$  
\cite{Are}.  The second Arens transpose of $m$   is the bounded bilinear map $$m^2: Y \times Z^* \ra X^*: (y, z^*) \mapsto y \c z^*$$ where $$\l m^2(y,z^*), x \r = \l z^*, m(x,y) \r; \ \ {\rm or} \ \ \l y \c z^*, x \r = \l z^*, x \c y \r. $$  More information about $m^2$,  and its relation with $m'$,  can be found in \cite{Sto1}.     \medskip 

We will often use the following readily verified  facts. The short argument found on page 309 of \cite{Eff-Rua} can be used to prove the first statement and the second statement is established in the proof of \cite[Proposition 7.1.2]{Eff-Rua}. 
 
\blem  \label{Transpose CC} Let $X, \ Y$ and $Z$ be operator spaces, $m: X \times Y \ra Z$ a completely contractive bilinear map. Then the following statements hold: \bi \item[(a)]  $m'$ and $m^2$ are  completely contractive; 
\item[(b)]  for  $x \in X$, the map $m_x: Y \ra Z$ defined by $m_x(y) = m(x,y)$ is c.b. with $\|m_x\|_{cb} \leq \|x\|$. \ei  \elem

 Suppose that $A$ is a Banach algebra with a contractive right approximate identity and $X$ is a contractive right Banach $A$-module via $(x,a) \mapsto x \c a$. Then the Cohen Factorization Theorem, as stated in Theorem 32.22 of \cite{HRII}, implies that $Z=X \c A$ is a closed $A$-submodule of $X$ and, moreover, 
 $$Z_{\| \c \| < 1} \subseteq X_{\| \c \| \leq 1} \c A_{\| \c \|\leq 1} \subseteq  Z_{\| \c \|\leq 1}.$$ 
(Taking $z \in Z$ with $\| z \| < 1$, $\delta = 1 - \| z \| >0$.  Then  $z = x \c a$ where $\| x - z \| \leq \delta $ -- hence $\| x\| \leq 1$ -- and $\| a \| \leq d =1$.)   We now state an operator space version of this result  due to P.J. Cohen on the factorization of modules. 

\blem \label{Operator Cohen Factorization Thm}  Suppose that $A$ is  a c.c. Banach algebra with a  contractive right approximate identity and $X$ is a c.c. right operator $A$-module via $(x,a) \mapsto x \c a$. Then $Z = X \c A$ is a closed $A$-submodule of $X$ and for each positive integer $n $, 
$$M_n(Z)_{\| \c \| < 1} \subseteq \left\{[x_{k,l} \c a] :  [x_{k,l}] \in M_n(X)_{\| \c \| \leq 1}, \   a \in  A_{\| \c \|\leq 1} \right\}\subseteq  M_n(Z)_{\| \c \|\leq 1}.$$ 

\elem

  \begin{proof}  As in the preceding paragraph,  $Z = X \c A$ is a closed $A$-submodule of $X$. Since $(x, a) \mapsto x \c a$ is a c.c., $([x_{k,l}], a ) \mapsto    [x_{k,l}] \c a:=  [x_{k,l} \c a] $ defines a contractive $A$-module action on $M_n(X)$ so,  by observing that $M_n(X) \c A = M_n(X \c A)$,  one sees that 
   the corollary follows from the phrasing of the Cohen Factorization Theorem provided above. 
  \end{proof}

When working in the category of operator spaces, we will often use the abbreviations c.b. in place of completely bounded and c.c. in place of completely contractive/complete contraction.  All undefined concepts from the theory of  operator spaces  can be found in  \cite{Eff-Rua} and \cite{Pis}.

\section{Arens product algebras and introverted homogeneous spaces} 

Let $A$ be a (completely contractive) Banach algebra. The following definition is motivated by the definition found on page 105 of \cite{Spr-Sto}. 

\bd   \label{Introverted Homogeneous Def}   \rm  A pair $(\SA, \| \c \|_\S)$, where $\SA$ is a linear subspace of $A^*$ and $\| \c \|_\S$ is a complete (operator space matrix) norm on $\SA$, will be called an \it  (operator) left homogeneous \rm subspace   of $A^*$  when the following two conditions are satisfied: 
\bi \item[(i)] the embedding $\SA \hookrightarrow A^*$ is a (complete) contraction;
\item[(ii)] $\SA$ is a right $A$-submodule of $A^*$ such that $\SA$ is a (completely) contractive right (operator) $A$-module with respect to $\| \c \|_\S$, i.e., 
\beq \label{Homogeneous Eqn}  \SA \times A \ra \SA \eeq 
is a (complete) contraction. 
\ei  
Observe that the first Arens transpose map of (\ref{Homogeneous Eqn}), 
\beq \label{Homogeneous Eqn Transpose}  \SA^*  \times \SA  \ra A^*: (\mu, \phi)  \mapsto \mu \sq \phi, \qquad  \l \mu \sq \phi, a \r = \l \mu , \phi \c a\r  \eeq
is also (completely) contractive by Lemma \ref{Transpose CC} .  We will say that $(\SA, \| \c \|_\S)$ is \it  (operator) left introverted \rm if  \bi \item[(iii)] the range of  the map described in  (\ref{Homogeneous Eqn Transpose}) is contained in $\SA$ and, moreover, 
\beq \label{Introverted Eqn}   \SA^*  \times \SA  \ra \SA : (\mu, \phi)  \mapsto \mu \sq \phi \eeq
is (completely) contractive with respect to $\| \c \|_\S$ and the associated dual operator space structure  -- which we will also  denote by $\| \c \|_\S$ --  on $\SA^*$.   \ei
(Operator) right introverted homogeneous subspaces of $A^*$ are similarly defined.  \ed 

 Note that a homogeneous subspace of $A^*$ will often fail to be closed in $A^*$.  When $\| \c \|_\S$ is the dual (operator space matrix) subspace norm $\|\c \|_{A^*}$ on  a closed subspace $\SA$ of $A^*$, our definition of a left introverted homogeneous space agrees with the the usual definition of a left introverted subspace of $A^*$ \cite{Dal-Lau, Fil-Mon-Neu}. The statement of the next proposition includes the introduction of some notation and terminology.

   
 \bp \label{Introverted Arens product algebras Prop} Let $(\SA, \| \c \|_\S)$ be a (operator) left introverted homogeneous subspace of $A^*$. 
 \bi \item[(a)] $(\SA^*, \| \c \|_\S, \sq)$  is a (c.c.)  Banach algebra with respect to the {\rm first Arens product}  $$\sq: \SA^* \times \SA^* \ra \SA^*: (\mu, \nu)\mapsto \mu \sq \nu \quad  defined \ by \quad \l \mu \sq \nu, \phi \r = \l \mu, \nu \sq \phi \r,$$
 and the map $$\eta_\S : A \ra \SA^* \quad  defined \ by \quad \l \eta_\S(a), \phi \r_{\S^*-\S} = \l \phi, a\r_{A^*-A}$$ is a (c.c.)  homomorphism with weak$^*$-dense range in $\SA^*$. Moreover, for each $\nu \in \SA^*$ and $a \in A$, the maps   
 $$ \mu \mapsto \mu \sq \nu \quad {\rm and} \quad  \mu \mapsto \eta_\S(a) \sq \mu$$ are  wk$^*$-wk$^*$ continuous on $\SA^*$. 
 \item[(b)] Let $\E_\S(A^*)$ denote the $\| \c \|_{A^*}$-closure of $\SA$ in $A^*$. Then $(\E_\S(A^*), \| \c \|_{A^*})$ is a (operator) left  introverted subspace of $A^*$ and 
 the embedding $\SA \hookrightarrow \E_\S(A^*)$ is a (complete) contraction.  
 \item[(c)]  Let $(\TA, \|\c \|_\T) $ be another (operator) left introverted homogeneous subspace of $A^*$. Then $\SA \subseteq \TA$ as a (complete) contraction if and only if there is a wk$^*$-continuous (complete) contraction $\Phi: \TA^* \ra \SA^*$ such that $\Phi \circ \eta_\T = \eta_\S$; the operator $\Phi$ is necessarily a wk$^*$-dense range homomorphism. 
 \ei 
 
 \ep

\begin{proof} We will provide the operator space version of the proof.

  (a)  Being the first Arens transpose of the map (\ref{Introverted Eqn}),  $\sq$  is  c.c. and bilinear. Associativity is readily established from the definition (and follows the same calculations used to establish associativity of $\sq$ on $A^{**}$ found, for example, in Section 2.6 of \cite{Dal}).   Observe that $\eta_\S = \iota^* \circ \, \widehat{}$ where  \ $\widehat{}: A \hookrightarrow A^{**}$ is the canonical embedding and $\iota: \SA \hookrightarrow A^*$. Since both  $\ \widehat{} \ $  and $\Phi = \iota^*$ have wk$^*$-dense range and $\Phi$ is wk$^*$-continuous, $\eta_\S$ also has wk$^*$-dense range and, as a composition of complete contractions, $\eta_\S$ is  c.c.. It is easy to directly check that $\eta_\S$ is a homomorphism (or note that both  $\ \widehat{} \ $ and --  using (c) --   $\Phi$ are homomorphisms).  

(b)  It is obvious that $\E_\S(A^*)$, the $\|\cdot\|_{A^*}$-closure of the $A$-submodule $\SA$ of $A^*$, satisfies both properties of an operator left homogeneous subspace of $A^*$.  It follows that $\sq: \E_\S(A^*)^* \times \E_\S(A^*) \ra A^*$ is  c.c., so we only need to establish that $\mu \sq \phi \in \E_\S(A^*)$ whenever $\mu \in \E_\S(A^*)^*$  and $ \phi \in \E_\S(A^*)$.  To see this, let $(\phi_n)$ be a sequence in $\SA$ such that $\| \phi_n - \phi \|_{A^*} \ra 0$. Then $\| \mu \sq  \phi_n - \mu \sq \phi \|_{A^*} \ra 0$ as well. Since $\iota: \SA \hookrightarrow \E_\S(A^*)$ is a (complete) contraction, $\iota^*(\mu) \in \SA^*$ and, from our assumption that $\SA$ is left introverted, $\mu \sq \phi_n = \iota^*(\mu) \sq \phi_n \in \SA$ for each $n$. Hence, $\mu \sq \phi \in \E_\S(A^*)$, as needed.   

(c) If the identity embedding $\iota: (\SA, \| \c \|_\S)\hookrightarrow (\TA, \| \c \|_\T)$ is a c.c., then it is easy to check that  its dual map $\Phi = \iota^*$ has all of the desired properties. Conversely, if $\Phi: \TA^* \ra \SA^*$ is a wk$^*$-continuous complete contraction  such that $\Phi \circ \eta_\T = \eta_\S$, then its predual map $\Phi_*: \SA \hookrightarrow \TA$ is the identity embedding.   \end{proof} 

When $(\SA^*, \sq)$ is a (operator) left introverted homogeneous subspace of $A^*$, we will refer to $(\SA^*, \sq)$ as a \it (operator) homogeneous left Arens product algebra over $A$\rm; as in \cite{Sto1}, we will call $(\SA^*, \sq)$ a \it (operator) left Arens product algebra over $A$ \rm when $\| \c \|_\S = \| \c \|_{A^*}$.     

\bd  \label{HLDBA Def}  \rm Let $\A$ be a (c.c.) Banach algebra with a fixed (operator space) predual $\A_*$.  The pair $(\A, \eta_\A)$ will be called a (operator) homogeneous left dual Banach algebra  -- (operator) HLDBA -- over $A$ if \bi \item[(i)]  for each $\nu \in \A$, $\mu \mapsto \mu \nu : \A \ra \A$  is  $\w*$ continuous; and
\item[(ii)] $\eta_\A: A \ra Z_t( \A)$  is a (completely)  contractive homomorphism  with weak$^*$-dense range in $\A$,  where $$Z_t(\A) = \{ \mu \in \A : \nu \mapsto \mu  \nu  \  {\rm is \ } \w* \ {\rm continuous \ on } \ \A\}$$ is the {\it topological centre} of $\A$.  \ei 
If, further, $\eta_\A(A_{\| \c \| \leq 1})$ is weak$^*$-dense in $\A_{\| \c \| \leq 1}$, then $(\A, \eta_A)$ is called a LDBA over $A$ \cite{Sto1}.\ed

Observe  that $\A_*$ is a left $\A$-submodule of $\A^*$, $Z_t(\A)$ is a norm-closed subalgebra of $\A$ and, viewing $\A$ as a  left $Z_t(\A)$-module via multiplication, $\A_*$ is a closed right $Z_t(\A)$-submodule of the  dual module $\A^*$.   If $(\A, \etaa)$ and $(\B,   \etab)$ are (operator) HLDBAs over $A$, we will write $(\A, \etaa) \geq (\B, \etab)$   if there is   a $\w*$ continuous (complete) contraction $\Phi: \A \ra \B$ such that $\Phi \circ \etaa = \etab$, and we call $\Phi$ a homomorphism of (operator) HLDBAs over $A$.  We call  $\Phi$ an isomorphism of (operator) HLBDAs over $A$ when it can be chosen to be a surjective (completely)  isometric isomorphism; in this case we   write     $(\A, \etaa) \cong (\B, \etab)$  and say that  $(\A, \etaa)$  and $(\B, \etab)$ are equivalent. 

 Observe that on any set of (operator) HLDBAs over $A$,  $\cong$ is an equivalence relation, $\leq$ is transitive,  our definition of $(\A, \etaa) \geq (\B, \etab)$ is consistent with  \cite[Def. 3.2]{Sto1}, and the intertwining map $\Phi$ is, as in the LDBA situation, necessarily a weak$^*$-dense range homomorphism.   However,  when $(\B, \etab)$ is not a LDBA over $A$, $\Phi$ may fail to be a surjection. 
 
  To see this, let  $G$ be an infinite locally compact group, $\lambda_{L^1(G)}: L^1(G) \ra B(L^2(G)) $ the left regular representation of $L^1(G)$, $VN(G)$ the von Neumann subalgebra of $B(L^2(G))$ generated by $\lambda_{L^1(G)}$. Then $(VN(G), \lambda_{L^1(G)})$ is an example of an HLDBA  over $L^1(G)$ that is not an LDBA (over $L^1(G)$). Since $C_0(G)$, the continuous functions on $G$ that vanish at infinity, is a closed introverted subspace of $L^\infty(G) = L^1(G)^*$, $(M(G), \eta_{C_0})$ is an LDBA over $L^1(G)$. Letting $\iota$ denote the identity embedding of the Fourier algebra $ A(G)$ into   $C_0(G)$, $\Phi = \iota^*: M(G) \ra VN(G)$ is a weak$^*$-continuous contraction such that $\Phi \circ \eta_{C_0} = \lambda_{L^1(G)}$; hence $(M(G), \eta_{C_0}) \geq (VN(G), \lambda_{L^1(G})$. Note however that $\Phi$    is  the left regular representation $\lambda_{M(G)}$  of $M(G)$, which is not surjective. Examples of this type are examined in greater generality in Sections 4 and 5. 

The Gelfand representation theorem shows that every right topological semigroup compactification of a locally compact group is a Gelfand compactication and, analogously,   \cite[Theorem 3.3]{Sto1} showed that every LDBA over $A$ is a left Arens product algebra over $A$ (and conversely), thus providing an abstract characterization of the left Arens product algebras.  
By Proposition \ref{Introverted Arens product algebras Prop},  every (operator) homogeneous left Arens product algebra over $A$ is a (operator) HLDBA over $A$.   We now establish the converse.

  \bt \label{HLDBA Rep Thm}   Every (operator) HLDBA over $A$ is equivalent to a unique (operator) homogeneous left Arens product algebra over $A$.    \et 

\begin{proof}  Let $(\A, \eta_\A)$ be an HLDBA over $A$. Let $\SA $ be the linear subspace $\etaa^* (\widehat{\A_*})$ of $A^*$ and consider the surjection $\Phi_*:= \eta_\A^* \circ \, \widehat{}: \A_* \ra \SA$, where $\, \widehat{}\,: \A_* \ra \A^*$ is the canonical embedding. Observe that $\Phi_*$ is also injective because $\eta_\A$ has weak$^*$-dense image in $\A$ and $\A$ separates points in $\A_*$. We can therefore define a complete (operator space matrix) norm $\| \c \|_\S$ on $\SA$ so that $\Phi_*$ is a (complete) isometry. The map $\Phi := (\Phi_*)^*: \SA^* \ra \A$ is then a weak$^*$-continuous (complete) isometry  such that $\Phi \circ \eta_\S = \etaa$.  In the operator space situation, we now show that $(\SA, \| \c \|_\S)$ satisfies the three axioms of Definition \ref{Introverted Homogeneous Def}. 

\smallskip 

\noindent (i) Let $\phi = (\Phi_*)_n(\psi) \in M_n(\SA)$, where $\psi = [\psi_{i,j}] \in M_n(\A_*)$. Observe that $\Phi_* = \eta_\A^* \circ \, \widehat{} \,$ is a complete contraction when viewed as a map of $\A_*$ into $A^*$, so 
$$\|\phi\|_{A^*} = \|(\Phi_*)_n(\psi) \|_{A^*} \leq \|\psi \|_{\A_*} = \| \phi \|_\S.$$ 

\smallskip 

\noindent (ii) For $a \in A$ and $\psi \in \A_*$, $\psi \c \eta_\A(a) \in \A_*$ since $\A_*$ is a right $Z_t(\A)$-submodule of $\A^*$. Moreover, a calculation shows that  $\Phi_*(\psi) \c a = \Phi_*(\psi \c \eta_\A(a))$, so $\SA$ is a right $A$-submodule of $A^*$.  To see that this module action is a c.c. with respect to $\| \c \|_\S$, let $\phi =  (\Phi_*)_r(\psi) \in M_r(\SA)$, where $\psi = [\psi_{i,j}] \in M_r(\A_*)$, and $a = [a_{k,l}] \in M_r(A)$. Since $\Phi_*$ is a complete isometry, $\A^*$ is a c.c. dual operator $\A$-module and $\eta_\A$  is a c.c., we obtain 
\beqs \|[\phi_{i,j} \c a_{k,l} ] \|_{\S} & = & \| [\Phi_*(\psi_{i,j} \c \eta_\A(a_{k,l}))]\|_{\S} = \| [\psi_{i,j} \c \eta_\A(a_{k,l})]\|_{\A_*}\\
& \leq &  \| [\psi_{i,j}]\|_{\A_*}  \|[ \eta_\A(a_{k,l})]\|_{\A} \leq   \| \psi\|_{\A_*}   \|[a_{k,l}]\|_A \\
& = & \| \phi \|_{\S} \| a \|_A
\eeqs
as needed. 

\smallskip 

\noindent (iii) For $\mu \in \SA^*$ and $\psi \in \A_*$, $\Phi(\mu) \c \psi \in \A_*$ and the calculation from the second paragraph of the proof of \cite[Thm. 3.3]{Sto1} shows that  $\mu \sq \Phi_*(\psi) = \Phi_*(\Phi(\mu) \c \psi)$ which belongs to $\SA$. Since $\Phi_*$ and $\Phi$ are complete isometries, and the left $\A$-module action on $\A_*$ is a c.c., an argument similar to the one used above to establish condition (ii)  shows that 
$$ \SA^*  \times \SA  \ra \SA : (\mu, \phi)  \mapsto \mu \sq \phi$$ 
is a complete contraction. 

\smallskip 

Hence, $(\SA^*, \eta_\S)$ is an operator homogeneous left Arens product algebra over $A$ and $(\A, \etaa) \cong (\SA^*, \eta_\S)$.  To establish uniqueness, suppose that $(\TA, \| \c \|_\T)$ is any operator left introverted homogeneous subspace of $A^*$ such that $ (\SA^*, \eta_\S) \cong (\TA^*, \eta_\T)$. Let $\Psi: \SA^* \ra \TA^*$  be a weak$^*$-continuous complete isometry such that $\Psi \circ \eta_\S = \eta_\T$.  Then, as noted in the proof of  Proposition \ref{Introverted Arens product algebras Prop}(c),  the predual map $\Psi_* : \TA \ra \SA$ is the identity embedding. Hence, $\TA = \SA$ and, since $\Psi_*$ is a complete isometry, the matrix norms $\| \c \|_\T$ and $\| \c \|_\S$ are equal. 
\end{proof} 

Given a (operator) HLDBA $(\A, \eta_\A)$ over $A$, we let $(\S_\A(A^*), \| \c \|_{\S_\A})$ be the (operator) left introverted homogeneous subspace of $A^*$ such that $(\A, \eta_A) \cong (\S_\A(A^*)^*, \eta_{\S_\A})$. Observe that by Theorem \ref{HLDBA Rep Thm} we can now, up to equivalence, view the class of all (operator) HLDBAs over $A$ as a set, denoted ${\cal HLD}(A)$, and $({\cal HLD}(A), \leq)$ is a partially ordered set.


\bc \label{HLDBA Ordering Cor}  Let $(\A, \eta_\A)$ and $(\B, \eta_\B)$ be (operator) HLDBAs over $A$. Then $(\B, \etab) \leq (\A, \etaa)$ if and only if $\S_\B (A^*)\subseteq \S_\A (A^*)$ and the embedding $\S_\B(A^*) \hookrightarrow \S_\A(A^*)$ is a (complete) contraction.  \ec 

\begin{proof} Since $(\B, \etab) \leq (\A, \etaa)$ exactly when $(\S_\B(A^*)^*, \eta_{\S_\B})   \leq (\S_\A(A^*)^*, \eta_{\S_\A})$, this is an immediate consequence of Proposition \ref{Introverted Arens product algebras Prop}(c). \end{proof}


\br  \label{Eberlein Min LDBA Remark} \rm Given any (operator) HLDBA $(\A, \eta_\A)$ over $A$, we observe  that $(\E_{\S_\A}(A^*)^*, \eta_{\E})$ is the minimum LDBA $(\B, \eta_\B)$  over $A$ such that $(\A, \eta_\A) \leq (\B, \eta_\B)$.  This is a consequence of Proposition \ref{Introverted Arens product algebras Prop}, Corollary \ref{HLDBA Ordering Cor} and \cite[Thm. 3.3]{Sto1}. \er 

The author has shown that $({\cal LD}(A), \leq)$ is a complete lattice by using the representation theorem  for LDBAs, \cite[Theorem 3.3]{Sto1}.    However, it is not obvious that Theorem \ref{HLDBA Rep Thm} can be used in a similar way to prove that 
$({\cal HLD}(A), \leq)$ is a complete lattice. With this in mind, and motivated by the corresponding notion found in the theory of   semigroup compactifications \cite{Ber-Jun-Mil},  we now introduce the construction  of the subdirect product  of a set of  (operator) HLDBAs.  Although we will focus our discussion on operator HLDBAs, it  will be clear that the construction also works in the category of Banach spaces.

Let  $\{ (\A_i, \eta_i): i \in I \}$ be a collection of operator HLDBAs over $A$. Let $ \prod_{i \in I} \A_i = \ell^\infty-\oplus_{i\in I} \A_i$ be the product operator space, as defined in \cite[Section 3.1]{Eff-Rua} and \cite[Section 2.6]{Pis}: 
$$ \|\mu \|_{M_n}= \sup_{i \in I} \|[\mu_{k,l}(i)]\|_{M_n(\A_i)}  \quad {\rm for} \quad \mu = [\mu_{k,l}] = [(\mu_{k,l}(i))_{i \in I} ] \in M_n\left(\prod_{i \in I} \A_i\right).$$  Then $ \prod_{i \in I} \A_i$ is the dual operator space of $\ell^1-\oplus (\A_i)_*$, e.g. see \cite[2.6.1]{Pis}, and with respect to the product defined by 
$$\mu \nu = (\mu_i)_{i \in I} (\nu_i)_{i \in I} = (\mu_i \nu_i)_{i \in I},$$ one can check that $\prod_{i \in I} \A_i$ is a c.c. Banach algebra. 

Define $\eta_\vee: A \ra \prod_{i \in I} \A_i$ by $\eta_\vee(a) = (\eta_i(a))_{i \in I}$, and let $\A_\vee$ denote the weak$^*$-closure of $\eta_\vee(A)$ in $\prod_{i \in I} \A_i$. Since each $\eta_i$ is a homomorphism, so is $\eta_\vee$.  For any $a \in M_n(A)$,   
\beqs \| (\eta_\vee)_n(a)\| & =&  \| [\eta_\vee(a_{k,l})] \| =  \| [(\eta_i(a_{k,l}))_{i \in I} ] \| \\ & =&    \sup_{i \in I} \| (\eta_i)_n(a)\|_{M_n(\A_i)} \leq  \| a\| \eeqs
since each $\eta_i$ is a complete contraction; hence, $\eta_\vee$ is a complete contraction. 

Observe that $$P_j : \prod_{i \in I} \A_i \ra \A_j: (\mu_i)_{i \in I} \mapsto \mu_j $$ is a weak$^*$-continuous completely contractive homomorphism, so 
$${\rm for \ each} \ j \in I \  \ P_j \mu_\alpha \ra P_j \mu \ \ {\rm weak}^* \ {\rm in} \ \A_j \quad {\rm whenever} \quad \mu_\alpha \ra \mu \ \ {\rm weak}^* \ {\rm in} \ \prod_{i \in I} \A_i.$$Moreover, it is easy to see that the converse of this last statement holds when the net $(\mu_\alpha)$ is bounded in $\prod_{i \in I} \A_i$.  Hence, a bounded mapping $\Phi: E^* \ra \prod_{i \in I} \A_i$ is weak$^*$-continuous on the unit ball of the dual Banach space $E^*$, and is therefore weak$^*$-continuous on all of $E^*$, if and only if $P_j \circ \Phi$ is weak$^*$-continuous on $E^*$ for each $j \in I$.  

Since each operator HLDBA $(\A_i, \eta_i)$ satisfies properties (i) and (ii) of Definition \ref{HLDBA Def}, it follows from these observations that for each $\nu \in \prod_{i \in I} \A_i$ and $a \in A$, the maps 
\beq \label{weak*-continuity of multn eqn}   
\mu \mapsto \mu \nu \qquad {\rm and} \qquad \nu \mapsto \eta_\vee(a) \nu
\eeq
are weak$^*$ continuous  on  $\prod_{i \in I} \A_i$. From this we see that $\A_\vee$ is a c.c. Banach subalgebra of $\prod_{i \in I} \A_i$: for $\mu, \nu \in \A_\vee$, if $\mu = {\rm w}^*-\lim \eta_\vee (a_\gamma) $ and $\nu = {\rm w}^*-\lim \eta_\vee (b_\delta) $, then for each $\gamma$, 
$$ \eta_\vee(a_\gamma) \nu = {\rm w}^*-\lim_\delta \eta_\vee (a_\gamma) \eta_\vee(b_\delta) =  {\rm w}^*-\lim_\delta \eta_\vee (a_\gamma b_\delta)  $$ belongs to $\A_\vee$, whence  $$\mu \nu = {\rm w}^*-\lim_\gamma \eta_\vee (a_\gamma)  \nu \in \A_\vee.$$ 
As a weak$^*$-closed subspace of the dual operator space $\prod_{i\in I} \A_i$, $\A_\vee$ with its subspace operator space structure is itself a dual operator space -- e.g., by \cite[Proposition 2.4.2]{Pis} -- and the weak$^*$-topology on $\A_\vee$ agrees with the relative weak$^*$-topology inherited from $\prod_{i \in I} \A_i$. Since the maps (\ref{weak*-continuity of multn eqn}) are weak$^*$-continuous on $\prod_{i \in I} \A_i$, they are also weak$^*$-continuous on $\A_\vee$. 

We have established that $(\A_\vee, \eta_\vee)$ is an operator HLDBA over $A$. We call $(\A_\vee, \eta_\vee)$ the \it subdirect product \rm of $\{(\A_i, \eta_i): i \in I\}$.  

\bt  \label{HLDBA Lattice Thm} The subdirect product $(\A_\vee, \eta_\vee)$ is the supremum of $\{(\A_i, \eta_i): i \in I\}$ in $({\cal HLD}(A), \leq)$. Hence, $({\cal HLD}(A), \leq)$ is a complete lattice with maximum element equal to $A^{**}$ and minimum element equal to the trivial Banach algebra.  \et 

\begin{proof}  Let $\Pi_j$ denote the restriction of $P_j$ to $\A_\vee$. Since $P_j$ is weak$^*$-continuous and a c.c., so is  $\Pi_j$. Moreover, $\Pi_j \circ \eta_\vee = \eta_j$, so $(\A_\vee, \eta_\vee) \geq (\A_j, \eta_j)$. Suppose now that $(\B, \eta_\B)$ is an operator HLDBA over $A$ such that for each $i \in I$, $(\B, \eta_\B) \geq (\A_i, \eta_i)$. Let $\Phi_i: \B \ra \A_i$ be a weak$^*$-continuous c.c. such that $\Phi_i \circ \eta_\B = \eta_i$ and define $$\Phi: \B \ra \prod_{i \in I} \A_i \quad {\rm by} \quad \Phi(\mu) = (\Phi_i(\mu))_{i \in I}.$$ Since each $\Phi_i$ is a c.c., the justification given above to show that $\eta_\vee$ is a c.c. shows that $\Phi$ is a c.c. as well. For each $j \in I$, $P_j \circ \Phi = \Phi_j$ is weak$^*$-continuous, which -- as noted above -- implies that $\Phi$ is weak$^*$-continuous on $\B$. Since $\Phi \circ \eta_\B = \eta_\vee$, weak$^*$-density of $\eta_\B(A)$ in $\B$ implies that $\Phi(\B)$ is contained in $\ov{\eta_\vee(A)}^{{\rm wk}^*} = \A_\vee$. Hence, $(\B, \eta_\B) \geq (\A_\vee, \eta_\vee)$.   
The infimum of $\{(\A_i, \eta_i): i \in I\}$ in $({\cal HLD}(A), \leq)$ is the supremum of the nonempty set of all lower bounds of $\{(\A_i, \eta_i): i \in I\}$.\end{proof} 

  The fundamental existence theorem for universal $P$-compactifications of locally compact groups \cite[Theorem 3.4]{Ber-Jun-Mil} demonstates the importance of subdirect products  in the construction of universal semigroup $P$-compactifications. We call an (operator) HLDBA $(\A, \eta_\A)$ over $A$ a $P$-extension of $A$ if $(\A, \eta_\A)$ has the property $P$ of (operator) HLDBAs, and say that $(\A, \eta_\A)$ is  \it a  \rm (or \it{the}, \rm up to equivalence)  \it universal $P$-extension of $A$ \rm if, further, $(\B, \eta_\B) \leq (\A, \eta_\A)$ whenever $(\B, \eta_\B)$ is a $P$-extension of $A$.  It seems worth noting that as  a corollary to Theorems \ref{HLDBA Rep Thm} and \ref{HLDBA Lattice Thm} we   have the following  version of  \cite[Theorem 3.4]{Ber-Jun-Mil} for (operator) HLDBAs. 

\bc Let $P$ be a property of (operator) HLDBAs over $A$ such $P$ is invariant under isomorphisms of (operator) HLDBAs. \bi
\item[(a)]  If $P$ is invariant under subdirect products, then $A$ has a universal $P$-extension. 
\item[(b)]  If $A$ has a universal $P$-extension and $P$ is invariant under homomorphisms of (operator) HLDBAs over $A$, then $P$ is invariant under subdirect products. 
 \ei
\ec 

\begin{proof} (a) Let $(\A_\vee, \eta_\vee)$ be the subdirect product of the set $\cal{P}$ of all (operator) homogeneous left Arens  product algebras $(\SA^*, \eta_\S)$ over $A$ for which $(\SA^*, \eta_\S)$ is a $P$-extension of $A$. By hypothesis $(\A_\vee, \eta_\vee)$  is a $P$-extension of $A$, and if $(\A, \eta_\A)$ is any $P$-extension of $A$, then by Theorem \ref{HLDBA Rep Thm}  $(\A, \eta_\A) $ is equivalent to some $(\SA^*, \eta_\S)$ in  $\cal{P}$.  By Theorem \ref{HLDBA Lattice Thm},  $(\A, \eta_\A) \leq (\A_\vee, \eta_\vee)$, so  $(\A_\vee, \eta_\vee)$ is the universal $P$-extension of $A$. \\
(b) Let $(\B, \eta_\B)$ be the universal $P$-extension of $A$ and let $(\A_\vee, \eta_\vee)$ be the subdirect product of a set $\{ (\A_i, \eta_i): i \in I\}$ of $P$-extensions of $A$. For each $i \in I$, $(\B, \eta_\B) \geq (\A_i, \eta_i)$, so $(\B, \eta_\B) \geq (\A_\vee, \eta_\vee)$ by Theorem \ref{HLDBA Lattice Thm}; by hypothesis $(\A_\vee, \eta_\vee)$ is a $P$-extension of $A$. \end{proof}

\section{The Fourier spaces $\FQA$ and $\FqA$}  
 
 Given any (completely contractive) Banach algebra $A$ and any (completely) contractive  right (operator) $A$-module action   
 $q: X \times A \ra X$ ($Q: X \times A \ra X$),   we will now introduce  the associated (operator) Fourier space $\FqA$  ($\FQA$). Since showing that $\FQA$ is an operator left introverted homogeneous subspace of $A^*$  involves more work than establishing  the corresponding statement for $\FqA$, we  will begin by focussing on the operator space situation.

\subsection{The operator Fourier space $\FQA$} 

Throughout this subsection, $A$ is a c.c. Banach algebra and $X$ is a c.c.  right operator $A$-module through the action $$Q: X \times A \ra X: (x,a) \mapsto x \c a.$$ By Lemma \ref{Transpose CC},  $$Q': X^* \times X \ra A^*: (x^*,x) \mapsto x^*\c x \quad  {\rm where} \quad   \l x^* \c x , a\r = \l x^*, x \c a\r$$
is also c.c.; thus, $Q'$ induces a complete contraction
 $$P: X^* \wot X \ra A^*: x^* \otimes x \mapsto x^* \c x.$$
Since $X$ is a c.c. right operator $A$-module and $X^*$ is a c.c. left operator (dual) $A$-module, the operator space projective tensor product $X^* \wot X$ becomes a c.c. operator $A$-bimodule in canonical fashion (cf. \cite[Theorem 2.6.4]{Dal}) and $P$ is an $A$-bimodule map.  Hence, the kernel $N$ of $P$ is a closed $A$-submodule of $X^* \wot X$ and the map $P_N$ 
determined by  the commuting diagram 
$$  \xymatrixrowsep{2pc} \xymatrixcolsep{3pc}
\xymatrix{ X^* \widehat{\otimes} X/ N  \ar@<.5ex>[rr]^{P_N}
 & & A^* 
  \\
X^* \widehat{\otimes} X  \ar@{->}[u]^{\Pi_N } \ar@<.5ex>[rru]_{P} } $$
is also a c.c. $A$-bimodule map into the dual $A$-bimodule $A^*$. Let $$\FQA := P(X^* \wot X) = P_N(X^* \wot X/N)$$
and give $\FQA$ the quotient operator space matrix norm, $\| \c \|_Q$, inherited from $X^* \wot X/N$ through the linear isomorphism $P_N$. Thus, for $\phi = [\phi_{i,j}] \in M_n(\FQA)$, 
\beq  \label{Q Norm Formula}     \| \phi \|_Q = \inf\{\| \xi \|_\wedge : \xi = [\xi_{i,j}] \in M_n( X^* \wot X) \ {\rm and} \ P_n(\xi) = \phi\}. \eeq 
We will call $\FQA$ the \it operator Fourier space \rm associated with $Q$ and $\| \c \|_Q$ is its \it Fourier operator space matrix norm.\rm 

\bt  \label{FQA operator introverted Thm} The pair $(\FQA, \| \c \|_Q)$ is an operator left introverted homogeneous subspace of $A^*$.    \et
 
 \begin{proof}  We will show that the three axioms of Definition \ref{Introverted Homogeneous Def}  are satisfied. 

 \medskip 
 
 \noindent (i)  Let $\phi = [\phi_{i,j}] \in M_n(\FQA)$, say $\phi = [P(\xi_{i,j})] = P_n(\xi) $.   Since $P: X^* \wot X \ra A^*$ is a c.c. map,  $\|\phi\|_{A^*} = \| P_n(\xi)]\|_{A^*} \leq \| \xi\|_\wedge$.  Hence,  $\| \phi \|_{A^*} \leq \| \phi\|_Q$ follows from the formula (\ref{Q Norm Formula}). 
 
 \medskip 
 
 \noindent (ii) Since $N$ is a closed $A$-submodule of $X^* \wot X$, it is easy to check -- and will be well known -- that $X^* \wot X/N$ is also a c.c. operator $A$-bimodule. As $P_N$  is a completely isometric $A$-bimodule isomorphism of $X^* \wot X/N$  onto $\FQA$, property (ii) is obvious. 
 
 \medskip 
 
 \noindent (iii) From (ii), $\FQA \times A \ra \FQA$ is a c.c. and therefore so is its first Arens transpose
 \beq   \label{FQA Introverted Eqn 1}   \sq_{A^*}: \FQA^* \times \FQA \ra A^*: (\mu, \phi) \mapsto \mu \sq \phi. \eeq 
 To establish (iii), we must first show that $\sq_{A^*}$ maps into $\FQA$.  To this end, observe that $Q': X^* \times X \ra \FQA$ is  c.c. since $P = P_N \circ \Pi_N: X^* \wot X \ra \FQA$ is a c.c. map. Hence
  \beq   \label{FQA Introverted Eqn 2}   \Psi = Q'': \FQA^* \times X^* \ra X^*: (\mu, x^*) \mapsto \mu \c x^* \ \ {\rm where} \ \ \l \mu \c x^*,x \r = \l \mu , x^* \c x \r \eeq
is also completely contractive. Fixing $\mu \in \FQA^*$, $\Psi_\mu(x^*) : = \Psi(\mu, x^*)$ defines a c.b. map on $X^*$ with $\|\Psi_\mu\|_{cb} \leq \| \mu\|$, whence $$\Psi_\mu \otimes {\rm id}_X : X^* \wot X \ra X^* \wot X$$ is also c.b..  We \it claim \rm that for any $\phi = P(\xi) \in \FQA$, $$\mu \sq \phi = P\circ (\Psi_\mu \otimes {\rm id}_X) (\xi) \in \FQA.$$ Since both $\xi \mapsto \mu \sq P(\xi) $ and $P \circ (\Psi_\mu \otimes {\rm id}_X)$ are continuous linear maps of $X^* \wot X$ into $A^*$ -- using (\ref{FQA Introverted Eqn 1}) and the continuity of $P$  when viewed as a map into either $\FQA$ or $A^*$ -- it suffices to establish the claim for $\xi = x^* \otimes x$.  To see this, observe that   
 \beq  \label{FQA Introverted Eqn 3}  \l \mu \sq (x^* \c x), a\r_{A^* -A} &  = &  \l \mu,   (x^* \c x) \c a\r_{\F_Q^* -\F_Q} =   \l \mu,   x^* \c (x \c a)\r_{\F_Q^* -\F_Q} \\
 & = & \l \mu \c   x^*,  x \c a\r_{X^*-X} = \l (\mu  \c x^*) \c x,  a\r_{A^*-A},  \nonumber \eeq  
 so \beqs  \mu \sq \phi = \mu \sq P(\xi) = \mu \sq (x^* \c x) = (\mu \c x^*) \c x = P (\Psi_\mu \otimes {\rm id}_X(\xi))\in \FQA 
 \eeqs
 as needed. 
 
 Finally, we will show that $$\sq: \FQA^* \times \FQA \ra \FQA$$ is c.c. with respect to $\| \c \|_Q$.  Let $$\widetilde{\Psi} : \FQA^* \wot X^* \ra X^* \quad {\rm and } \quad \widetilde{\sq}_{A^*} : \FQA^* \wot \FQA \ra A^*$$ be the c.c. linearization mappings of $\Psi$ and $\sq_{A^*}$ respectively. From (i), each of the maps in the diagram

 $$\xymatrixrowsep{2pc} \xymatrixcolsep{4pc}
\xymatrix{ \FQA^* \wot (X^* \wot X)  \ar@<.5ex>[dd]_{{\rm id} \otimes P} \ar@<.5ex>[r]^{T}_{\cong}   & (\FQA^* \wot X^*) \wot X   \ar@<.5ex>[r]^-{  \widetilde{\Psi} \otimes {\rm id}_X }   &  X^* \wot X  \ar@<.5ex>[d]^{P}  
  \\ && \FQA   \ar@{^{(}->}[d]^{ \iota}    \\
\FQA^* \wot \FQA  \ar@<.5ex>[rr]^{\widetilde{\sq}_A^*}  &&     A^*  }$$
is a c.c.  Moreover, it follows from (\ref{FQA Introverted Eqn 3}) that the diagram commutes on elementary tensors $\mu \otimes (x^* \otimes x)$ and therefore commutes on all of $\FQA^* \wot (X^* \wot X)$.  Thus, for $(\mu, \xi) \in \FQA^* \times (X^* \wot X)$, $\mu \sq P\xi \in \FQA$ and $$\mu \sq P\xi = \widetilde{\sq}_{A^*}\circ ( {\rm id} \otimes P) (\mu  \otimes \xi) = P \circ (\widetilde{\Psi} \otimes {\rm id}_X) \circ T (\mu \otimes \xi).$$ Let $r$ and  $s$ be positive integers, and consider 
$$  \sq_{r;s}: M_r(\FQA^*)  \times M_s(\FQA)  \ra M_{r \times s}(\FQA).$$  Let  $\mu = [\mu_{i,j}] \in M_r(\FQA^*)$,  $\phi = [\phi_{k,l}] \in M_s(\FQA)$, and take any   $\xi = [\xi_{k,l}] \in M_s(X^* \wot X)$  such that $P_s \xi = \phi$. Then
\beqs  \| \sq_{r;s} (\mu, \phi) \|_Q & = & \| [\mu_{i,j} \sq P\xi_{k,l} ]\|_Q = \| P_{r \times s} \circ (\widetilde{\Psi} \otimes {\rm id}_X)_{r \times s} \circ T_{r \times s} [\mu_{i,j} \otimes \xi_{k,l} ]\|_Q \\
& \leq & \| [\mu_{i,j} \otimes \xi_{k,l}]\|_\wedge = \| \mu \otimes \xi \|_\wedge = \|\mu \| \|\xi\|_\wedge.
\eeqs
 The formula (\ref{Q Norm Formula}) now yields $\|  \sq_{r;s} (\mu, \phi) \|_Q \leq \| \mu \| \| \phi\|_Q$. Hence, $\sq$ is a complete contraction. 
 \end{proof} 
 
 Hence,  the Fourier dual  $(\FQA^*, \sq)$ is an operator  homogeneous left Arens product algebra  over $A$. We will denote the embedding homomorphism $\eta_{\F_Q}: A \ra \FQA^*$ by $\eta_Q$.   We begin our study of $\FQA^*$  by showing that it can be identified with a weak$^*$-closed subalgebra of $CB(X^*)$, where the weak$^*$-topology on $CB(X^*)$ is, as usual, defined through its canonical identification with $(X^* \wot X)^*$ \cite[Cor. 7.1.5]{Eff-Rua}.   
 
 The left dual $A$-module action on $X^*$,
 $Q^2: A \times X^* \ra X^*$ (as defined in Section 2), is a c.c. map and so $$\Gamma_Q: A \ra CB(X^*), \quad \Gamma_Q(a)(x^*) = a \c x^*$$ defines a representation of $A$ on $X^*$ such that $\|\Gamma_Q(a) \|_{cb} \leq \| a\|$. Let 
 $$M_Q =  \overline{\Gamma_Q(A)}^{{\rm wk}^*}\preceq CB(X^*)$$ 
 be the weak$^*$-closed operator subalgebra of $CB(X^*)$  generated by $\Gamma_Q$.   As noted in the proof of Theorem \ref{FQA operator introverted Thm} -- see (\ref{FQA Introverted Eqn 2}) --  $$\Psi: \FQA^* \times X^* \ra X^*: (\mu, x^*) \mapsto \mu \c x^*$$ 
 is a completely contractive bilinear map. It follows that $\bg_Q$ maps $\FQA^*$ into $CB(X^*)$ where 
 $$\bg_Q(\mu) (x^*) = \mu \c x^*$$ 
 and $\|\bg_Q(\mu)\|_{cb} \leq \| \mu \|$.  We now observe that  $\bg_Q$ is a weak$^*$-continuous completely isometric extension of $\Gamma_Q$ to a representation of $\FQA^*$ as c.b. operators  on $X^*$.  
 
 \bt  \label{Completely Isometric Repn Thm}  The following statements hold.  
 \bi \item[(a)] The bilinear map 
 $$ \FQA^* \times X^* \ra X^*: (\mu, x^*) \mapsto \mu \c x^* \ \   where \ \ \l\mu \c x^*, x \r = \l \mu, x^* \c x\r$$ 
 defines a c.c. left operator $\FQA^*$-module action on $X^*$ such that for each $x^* \in X^*$ 
 $$\mu \mapsto \mu \c x^*: X^* \ra X^*$$ is weak$^*$-continuous and $\eta_Q(a) \c x^* = a \c x^*$. 
 \item[(b)]  The map $\bg_Q: \FQA^* \ra CB(X^*)$ is a weak$^*$-continuous completely isometric algebra isomorphism of $\FQA^*$ onto $M_Q$ such that $\bg_Q \circ \eta_Q = \Gamma_Q$. 
 \ei
   \et 
   
 \begin{proof}  (a) We have already observed that the bilinear map $(\mu, x^*) \mapsto \mu \c x^*$ is a complete contraction. For $\mu, \nu \in \FQA^*$, 
 $$\l (\mu \sq \nu) \c x^*, x \r = \l \mu \sq \nu, x^* \c x \r = \l \mu, \nu \sq (x^* \c x) \r = \l \mu, (\nu \c x^*) \c x \r = \l \mu \c (\nu \c x^*), x \r,$$ where we have used the calculation  (\ref{FQA Introverted Eqn 3}).  Also, 
 $$\l \eta_Q(a) \c x^*, x \r = \l \eta_Q(a), x^* \c x \r = \l x^* \c x, a\r = \l x^*, x \c a \r = \l a\c x^*, x \r.$$

(b) From our definition of $\| \c \|_Q$, $P : X^* \wot X \ra \FQA$ is a complete quotient map, so $P^* : \FQA^* \ra (X^* \wot X)^*$ is a weak$^*$-continuous complete isometry, as is the canonical mapping $\lambda : (X^* \wot X)^* \ra CB(X^*)$. For $\mu \in \FQA^*$, 
$$\l \lambda \circ P^*(\mu)(x^*), x \r = \l P^*(\mu), x^* \otimes x \r = \l \mu, x^* \c x \r = \l \mu \c x^*, x \r,$$
so $\bg_Q = \lambda \circ P^* $ is a weak$^*$-continuous complete isometry of $\FQA^*$ into $CB(X^*)$. From part (a), $\bg_Q$ is an algebra isomorphism and $\bg_Q \circ \eta_Q = \Gamma_Q$. 
Since $\bg_Q$ is weak$^*$-continuous with norm-closed range, its range is also weak$^*$-closed in $CB(X^*)$, see e.g. \cite[Theorem VI.1.10]{Con}. As $\Gamma_Q(A) = \bg_Q(\eta_Q(A))$, we can conclude that $M_Q$ is contained in  ${\rm range}(\bg_Q)$. By Proposition \ref{Introverted Arens product algebras Prop},  $\eta_Q$ has weak$^*$-dense range in $\FQA^*$. The containment of ${\rm range}(\bg_Q)$ in $ M_Q$ follows because $\bg_Q$ is weak$^*$-continuous and $\bg_Q \circ \eta_Q = \Gamma_Q$.     
 \end{proof} 
 
 \bc \label{FQA*=MQ Identification Cor}  We can identify the c.c. Banach algebra $(\FQA^*, \sq)$ with the weak$^*$-closed operator subalgebra  $M_Q$ of $CB(X^*)$ via a weak$^*$-homeomorphic completely isometric algebra isomorphism 
 $$M_Q \ra \FQA^*: T \mapsto \mu_T \quad   satisfying   \quad \l \mu_T, x^* \c x \r = \l Tx^*, x\r.$$
 \ec

 \bp \label{M_Q Description Prop}  Let   $T \in CB(X^*)$ and consider the following statements: 
 \bi 
 \item[(a)] $T \in M_Q$;
 \item[(b)] there is a c.b. linear map $\Lambda_T: \FQA  \ra \FQA$ such that $\Lambda_T(x^* \c x) = (Tx^*) \c x$ and $\|\Lambda_T \|_{cb} \leq \|T\|_{cb}$; 
   \item[(c)] there is a bounded linear map $\Lambda_T: \FQA  \ra \FQA$ such that $\Lambda_T(x^* \c x) = (Tx^*) \c x$.  
 \ei 
 Then (a) implies (b)  implies (c); and (c) implies (a) whenever $A$ has a  bounded approximate identity for $X$ (e.g., when $A$ has a right bounded approximate identity  and $X$ is essential). Moreover, when $T \in M_Q$, $\Lambda_T (\phi) = \mu_T \sq \phi$ and $\Lambda_T = (R_T)_*$ where $R_T: M_Q \ra M_Q: S \mapsto ST$.  
  \ep 
  
  \begin{proof}  Suppose that $T \in M_Q$.  By Lemma \ref{Transpose CC}, the map $$\Lambda_T: \FQA \ra \FQA: \phi \mapsto \mu_T \sq \phi$$ is c.b. with $\|\Lambda_T \|_{cb} \leq  \| \mu_T\| = \|T\|_{cb}$. Observe that 
  \beqs \l \Lambda_T(x^* \c x), a \r & = & \l \mu_T \sq (x^* \c x), a\r = \l \mu_T, (x^* \c x) \c a\r  = \l \mu_T, x^* \c (x \c a) \r \\
  & = & \l Tx^*, x \c a \r = \l (Tx^*) \c x , a \r 
  \eeqs
  so $\Lambda_T(x^* \c x) = (Tx^*) \c x$.  Also, $\Lambda_T^*: M_Q \ra M_Q$ and for $S \in M_Q$, $\phi \in \FQA$, 
  $$ \l \Lambda_T^*(S), \phi \r = \l\mu_S, \mu_T \sq \phi \r = \l \mu_S \sq \mu_T, \phi \r = \l \mu_{ST}, \phi \r = \l ST, \phi \r,$$ so $\Lambda_T^*(S) = ST = R_T(S)$.  This establishes the first implication and $(b)$ implies $(c)$ is trivial. 
  
  Assuming that statement $(c)$ holds and $(e_i)$ is a bounded approximate identity for the right $A$-module $X$, let $\mu$ be a weak$^*$-limit point of the net $(\Lambda_T^*(\eta_Q(e_i)))$ in $\FQA^*$; we can assume  that   $\Lambda_T^*(\eta_Q(e_i)) \ra \mu$ weak$^*$ in $\FQA^*$. Then 
 \beqs \l \bg_Q(\mu)(x^*), x \r & = & \l \mu \c x^*, x \r = \l \mu, x^* \c x \r \\
 & = & \lim \l \Lambda_T^*(\eta_Q(e_i)), x^* \c x \r = \lim \l \eta_Q(e_i), (Tx^*) \c x \r \\
 & = & \lim \l (Tx^*)  \c x, e_i \r = \lim \l Tx^*, x \c e_i \r \\
 &=& \l Tx^*, x \r. 
 \eeqs 
 Hence, $T = \bg_Q(\mu)$ belongs to $ M_Q$.   
  \end{proof}
  
  Let $\EQA$ denote the $\|\c\|_{A^*}$-closure of $\FQA$ in $A^*$, i.e., let $\EQA= \E_{\F_Q}(A^*)$. Then by   Remark  \ref{Eberlein Min LDBA Remark}, the \it Eberlein space \rm associated with $Q$, $(\EQA, \| \c \|_{A*})$ is an operator left introverted subspace of $A^*$ such that  $(\FQA^*, \eta_Q) \leq (\EQA^*, \eta_\E)$. Moreover,  $(\EQA^*, \eta_\E)$ is the minimum LDBA over $A$ with this property.
 
 \br \label{FQA and EQA Span Remark}   \rm (a)  Observe  that the operator Fourier space $\FQA$ is the $\|\c\|_Q$-closed linear span of $$X^* \c X = \{ x^* \c x : x^* \in X^*, \ x \in X\}$$ and $\EQA$ is the $\| \c \|_{A^*}$-closed linear span of $X^* \c X$. Indeed, since $P: X^* \wot X \ra \FQA$ is $\| \c \|_Q$-continuous and  ${\rm span}\{X^* \c X\} = P (X^* \otimes X)$, this easily follows from the definitions.
 \smallskip 
 \noindent (b) When $X$ is an essential c.c. right operator $A$-module,  $\FQA$ and $\EQA$ are contained in $LUC(A^*)$, the $\|\c\|_{A^*}$-closed linear span of $A^* \c A$.  Indeed, as noted in \cite[Example 4.3(b)]{Sto1},  $x^* \c (x \c a) = (x^* \c x) \c a \in A^* \c A $, so this follows from (a).  In particular, suppose that $G$ is a locally compact group and $X$ is a right Banach $G$-module, and $X$ is therefore a neo-unital right Banach $L^1(G)$-module through the action defined by the weak integral $x\c f = \int x \c s \, f(s) \, ds$. Assuming that this $L^1(G)$-module action is c.c.,  ${\cal F}_Q(L^1(G)^*)$ and ${\cal E}_Q(L^1(G)^*)$ are contained in $LUC(G)$ and, as observed in Section 4.3 of \cite{Sto1}, in this case $x^* \c x \in LUC(G)$ is  given by 
 $$(x^* \c x)(s) = \l x^*, x \c s \r \quad (s \in G).$$  
 
 \noindent (c) Suppose that $(\TA, \| \c \|_\T)$ is an operator left introverted homogeneous subspace of $A^*$ such that $(\FQA^*, \eta_Q) \leq (\TA^*, \eta_\T)$. Then  
 \beq   \label{Repn of TA* Eqn} \Gamma_\T : \TA^* \ra CB(X^*) \quad {\rm given \ by} \quad \l \Gamma_\T(\mu)(x^*), x \r = \l \mu, x^* \c x \r \eeq
 defines a c.c. weak$^*$-continuous representation of $\TA^*$ on $X^*$ such that $\Gamma_\T\circ\eta_\T(a)(x^*) = a \c x^*$. To see this, observe that $\Phi: \TA^* \ra \FQA^*$, the dual of the embedding map $\FQA \hookrightarrow \TA$, is a weak$^*$-continuous c.c. homomorphism such that $\Phi\circ \eta_\T = \eta_Q$ and $\bg_Q \circ \Phi = \Gamma_\T$ (see Corollary \ref{HLDBA Ordering Cor} and Proposition \ref{Introverted Arens product algebras Prop}(c)).  Since $\FQA \hookrightarrow \EQA$ has dense range, $\Phi: \EQA^* \ra \FQA^*$ is a monomorphism, and therefore $\bg_Q \circ \Phi$ is one-to-one. Thus, 
 \bi \item[]  (\ref{Repn of TA* Eqn})  determines a \it faithful \rm  c.c. weak$^*$-continuous representation of  $\EQA^*$ on $X^*$. 
 \ei  
 If $X$ is an essential right Banach $A$-module through $Q$ then, as noted in (b), $(\FQA^*, \eta_Q) \leq (LUC(A^*)^*, \eta_{LUC})$, so   \bi \item[]  (\ref{Repn of TA* Eqn})  determines a c.c. weak$^*$-continuous representation, $\Gamma_{LUC}$,  of  $LUC(A^*)^*$ on $X^*$ in this case.  
 \ei   
 Moreover, if $X$ is an essential $L^1(G)$-module and $\Theta: M(G) \hookrightarrow LUC(G)^*$ is the canonical  embedding given by $\l \Theta(\mu), f \r = \int f \, d\mu$, then $(\Gamma_{LUC} \circ \Theta) (\mu)(x^*) = \mu \c x^*$, the usual left dual module  action of $M(G)$ on $X^*$. That is, $\Gamma_{LUC}$ extends the usual left dual $M(G)$-module action on $X^*$ to a weak$^*$-continuous module action of $LUC(G)^*$ on $X^*$; a different justification of this statement can be found in  \cite{Sto1}.  
 \er

   \bt \label{EQA Min LDBA Prop}   The following statements hold. \bi \item[(a)] In $({\cal HLD}(A), \leq )$ [$({\cal LD}(A), \leq )$], $(\FQA^*, \eta_Q)$ [$(\EQA^*, \eta_\E)$] is the minimum operator HLDBA  [LDBA] $(\A, \eta_\A)$ over $A$ for which there is a c.c. weak$^*$-continuous representation $\Gamma_\A : \A \ra CB(X^*)$ extending $\Gamma_Q$ in the sense that $\Gamma_\A \circ \eta_\A = \Gamma_Q$. 
   \item[(b)] $(\FQA^*,\eta_Q)$ [$(\EQA^*, \eta_\E)$] is the unique operator HLDBA [LDBA] over $A$ for which there is a completely isometric [c.c. faithful] weak$^*$-continuous representation $\Gamma_\A: \A \ra CB(X^*)$ such that $\Gamma_\A \circ \eta_\A = \Gamma_Q$. 
   \ei   \et 
 
 \begin{proof}  (a) By Theorem \ref{Completely Isometric Repn Thm}, $(\FQA^*, \eta_Q)$ has this property.  Suppose that $(\A, \eta_\A)$ is an operator HLDBA over $A$ with this property. Assuming without loss of generality that $(\A, \eta_\A)$ is an operator homogeneous left Arens product algebra $(\SA^*, \eta_\S)$ over $A$ (Theorem \ref{HLDBA Rep Thm}), let $\Gamma_\S: \SA^* \ra CB(X^*)$ be a weak$^*$-continuous c.c. representation of $\SA^*$ such that $\Gamma_\S \circ \eta_\S = \Gamma_Q$.   Letting $\sigma: X^* \wot X \ra \SA$ be the c.c. predual map of $\Gamma_\S$, we have 
 \beqs  \l \sigma(x^* \otimes x), a\r_{A^*-A} & = & \l \eta_\S(a), \sigma(x^* \otimes x) \r = \l \Gamma_\S \circ \eta_\S(a), x^* \otimes x \r  \\
 & = &  \l \Gamma_Q(a)( x^*), x \r  = \l x^* \c x , a\r_{A^*-A}. 
 \eeqs
 Hence, $\sigma: X^* \wot X \rightarrow \SA: x^* \otimes x \mapsto x^* \c x$; since ${\rm id}_\S: \SA \hookrightarrow A^*$ is also a c.c. -- see Definition \ref{Introverted Homogeneous Def} -- we obtain  $P = {\rm id}_\S \circ \sigma: X^* \wot X \ra A^*$.   Hence, $\FQA = P(X^* \wot X) \subseteq \SA$,   $N = \ker P = \ker \sigma$, and we obtain a c.c. $\sigma_N : X^* \wot X /N \ra \SA$; thus, $\sigma_N \circ P_N^{-1}: \FQA \ra \SA$ is a c.c.. A calculation shows that  $\sigma_N \circ P_N^{-1}$ is the identity embedding $\FQA \hookrightarrow \SA$, so  $(\FQA^*, \eta_Q) \leq (\SA^*, \eta_\S)$ by Corollary \ref{HLDBA Ordering Cor}. 
 
 As noted in Remarks \ref{FQA and EQA Span Remark}(c), $(\EQA^*, \eta_\E)$ satisfies the desired property. If $(\A, \eta_\A) \cong (\SA^*, \eta_\S)$ is a LDBA over $A$ satisfying this property, then $(\FQA^*, \eta_Q) \leq (\SA^*, \eta_\S)$ by the  case above, and therefore $\FQA \subseteq \SA$ by Corollary \ref{HLDBA Ordering Cor}. Since $\SA$ is a closed subspace of $A^*$, $\EQA \subseteq \SA$ and it follows that $(\EQA^*, \eta_\E) \leq (\SA^*, \eta_\S)$, using \cite[Corollary 3.5]{Sto1}.
 \smallskip
 
 \noindent (b) By Theorem \ref{Completely Isometric Repn Thm},  $\bg_Q$ is a completely isometric representation of $\FQA^*$ mapping weak$^*$-homeomorphically onto $M_Q$, from which it follows  that $(M_Q, \Gamma_Q)$ is an operator HLDBA over $A$ and $(\FQA^*, \eta_Q) \cong (M_Q, \Gamma_Q)$.  If $(\A, \eta_\A)$ is any operator HLDBA over $A$ for which there exists a weak$^*$-continuous completely isometric representation 
  $\Gamma_\A: \A \ra CB(X^*)$ such that $\Gamma_\A \circ \eta_\A = \Gamma_Q$, then the argument provided in the last paragraph of the proof of Theorem \ref{Completely Isometric Repn Thm} shows that $\Gamma_\A$ maps onto $M_Q$. Hence $(\A, \eta_\A) \cong (M_Q, \Gamma_Q)$ as well, so $(\A, \eta_\A)\cong (\FQA^*, \eta_Q)$.  
  
  In Remarks \ref{FQA and EQA Span Remark}(c), we observed that  $(\EQA^*, \eta_\E)$ has a weak$^*$-continuous c.c. faithful representation $\Gamma_\E$ on $X^*$ such that $\Gamma_\E \circ \eta_\E = \Gamma_Q$. If $(\A, \eta_\A) \cong (\SA^*,  \eta_\S)$ is any LDBA over $A$ with a weak$^*$-continuous c.c. faithful representation $\Gamma_\S: \SA^* \ra CB(X^*)$ such that $\Gamma_\S \circ \eta_\S = \Gamma_Q$, then $(\EQA^*, \eta_\E)\leq(\SA^*,  \eta_\S)$ by (a). This means that $\EQA \subseteq \SA$ and, letting $\iota$ denote the associated embedding map, $\Phi = \iota^*: \SA^* \ra \EQA^*$  satisfies $\Phi \circ \eta_\S = \eta_\E$. Since $\Gamma_\E \circ \Phi \circ \eta_\S = \Gamma_\E \circ \eta_\E = \Gamma_Q = \Gamma_\S \circ \eta_\S$, weak$^*$-density of $\eta_\S(A)$ in $\SA^*$ implies that $\Gamma_\E \circ \Phi = \Gamma_\S$. Hence, $\Phi$ is injective and, as noted in Definition 3.2 of \cite{Sto1}, $\Phi$ is necessarily surjective; therefore $\iota$ maps onto $\SA$. Thus, $\EQA = \SA$, as needed.  \end{proof} 

Observe that for each $a \in A$,  $\Gamma_Q(a)(x^*) = a \c x^*$ is a weak$^*$-continuous c.b. operator on $X^*$, i.e., observe that $\bg_Q \circ \eta_Q = \Gamma_Q$ maps $A$ into   $CB^\sigma(X^*)$.

\bc  \label{HLDBA=FQ*Cor}  Let $(\A, \eta_\A)$ be an operator  HLDBA [LDBA] over $A$ for which there is a completely isometric [c.c. faithful] weak$^*$-continuous representation $\Gamma_\A: \A \ra CB(X^*)$ such that  $\Gamma_\A \circ \eta_\A$ maps $A$ into $CB^\sigma(X^*)$. Then there is a c.c. right operator $A$-module action $Q: X \times A \ra X$ such that $\Gamma_\A \circ \eta_\A = \Gamma_Q$ and $(\A, \eta_\A) \cong (\FQA^*, \eta_Q)$   [$\cong (\EQA^*, \eta_\E)$].   \ec 

\begin{proof}  Since $\Gamma_\A \circ \eta_\A \in CB(A, CB(X^*))$ is a c.c. homomorphism, $X^*$ is a c.c. left operator $A$-module through 
$$\gamma: A \times X^* \ra X^*: (a, x^*) \mapsto a \c x^*:= \Gamma_\A \circ \eta_\A(a)(x^*)$$
\cite[Prop. 7.12]{Eff-Rua}.  For each $a \in A$, let $$Q(\c, a): X \ra X : x \mapsto x \c a$$ denote the predual map of the weak$^*$-continuous map $\gamma(a, \c)$ on $X^*$.   Then $\l a \c x^*, x\r = \l x^*, x \c a\r$, from which it follows that  $ x \c (ab) = (x \c a) \c b$. Observe that $$\gamma^2: X^* \times X^{**} \mapsto A^*  \quad {\rm and} \quad \gamma^{22} : X^{**} \times A^{**} \ra X^{**}$$
are c.c. bilinear (e.g.,  by Lemma \ref{Transpose CC})  and $\gamma^{22} (\widehat{x}, \widehat{a}) = Q(x, a)^\wedge$.  It follows that $X$ is a c.c. right operator $A$-module through $Q$ such that $\gamma$ is the associated left dual module action on $X^*$. Since $\Gamma_\A \circ \eta_\A(a)(x^*) = a \c x^* = \Gamma_Q(a)(x^*)$,  $(\A, \eta_\A) \cong (\FQA^*, \eta_Q)$   [$\cong (\EQA^*, \eta_\E)$] by Theorem \ref{EQA Min LDBA Prop}.   \end{proof}

\br \label{FQA=EQARemark}  \rm It follows from Corollary \ref{HLDBA=FQ*Cor} that several commonly studied  operator HLDBAs  over $A$ can be recognized as   $(\FQA^*, \eta_Q)$ or $(\EQA^*, \eta_\E)$  for some module action $Q$. We discuss some specific examples in the next section. Moreover, if $(\A, \eta_\A)$ is an operator \it LDBA  \rm over $A$ for which there is a completely isometric weak$^*$-continuous representation $\Gamma_\A: \A \ra CB(X^*)$ such that $\Gamma_\A \circ \eta_\A(A)$ is contained in $CB^\sigma(X^*)$, then there is a c.c. right operator $A$-module action $Q$ on $X$ such that $\Gamma_\A \circ \eta_\A = \Gamma_Q$,  $(\A, \eta_\A) \cong (\EQA^*, \eta_\E) \cong (\FQA^*, \eta_Q)$ and -- by Corollary \ref{HLDBA Ordering Cor} -- $\FQA =\EQA = \S_\A(A^*)$ with equality of the matrix norms  $\| \c \|_Q$ and $ \| \c \|_{A^*}$  on this common space 
\er 
  
   It will be useful for us to be able to directly determine when $\FQA = \EQA$ and $\|\c \|_Q = \| \c \|_{A^*}$ without employing Corollary  \ref{HLDBA=FQ*Cor}  together with pre-existing theory. 
 In the following proposition,  $\iota$ denotes  the c.c. embedding $ \FQA \hookrightarrow \EQA $ and $\Phi:= \iota^*: \EQA^* \ra \FQA^*$ is the weak$^*$-continuous c.c. weak$^*$-dense range homomorphism such that $\Phi \circ \eta_\E = \eta_Q$.   
  
  \bp \label{FQA =EQA Prop} The following statements are equivalent: 
  \bi \item[(a)] $\FQA = \EQA$ and the matrix norms $\| \c \|_Q $ and $\| \c \|_{A^*}$ are equal on this common space;
  \item[(b)] $\Phi: \EQA^* \ra \FQA^*$ is a complete isometry; 
  \item[(c)] for each positive integer $n$ and each $\mu \in M_n(\EQA^*)$, $\|\Phi_n(\mu)\|_Q \geq \|\mu\|_{A^*}$. 
  \ei 
  These three equivalent conditions hold whenever 
  $${\cal S}_n = \left\{ [x_{i,j}^* \c x_{k,l}] : [x_{i,j}^*] \in M_r(X^*)_{\| \c \| \leq 1}, \ [x_{k,l}] \in M_s(X)_{\| \c \| \leq 1} \ {\rm and} \ rs = n\right\}$$ 
  is $\|\c \|_{A^*}$-dense in $M_n(\EQA)_{\| \c \| \leq 1}$ for each positive integer $n$.   \ep 
 
 \begin{proof}   The implication $(a)$ implies $(b)$ is obvious, as is the equivalence of statements $(b)$ and $(c)$ since $\Phi$ is a c.c.. Suppose that $\Phi = \iota^*$ is a complete isometry. Then $\iota$ has closed dense range in $\EQA$ and it follows that $\FQA= \EQA$. By the open mapping theorem, the norms $\| \c \|_Q$ and $\|\c \|_{A^*}$ are equivalent on this set, i.e. on $M_1(\EQA)$, so the sets $\FQA^*$ and $\EQA^*$ are also one and the same. Hence, $\Phi= \iota^*$ is the identity mapping of $\EQA^*$ onto $\FQA^*$. As $\Phi$ is a complete isometry, so is $\iota= \Phi^*|_{\FQA}$, which establishes the final implication $(b)$ implies $(a)$. 
 
 Assuming that ${\cal S}_n$ is $\| \c \|_{A^*}$-dense in $M_n(\EQA)_{\| \c \| \leq 1}$, we will now show that condition (c) holds. To see this, first note that ${\cal S}_n$ is contained in $M_n(\FQA)_{\| \c \|_Q \leq 1}$ since, as observed in the proof of Theorem \ref{FQA operator introverted Thm}, $Q': X^* \times X \ra \FQA$ is a c.c. with respect to $\| \c \|_Q$.  Taking $\mu = [\mu_{e,f}] \in M_n(\EQA^*)$,  this observation and our hypothesis give 
 \beqs \|\Phi_n(\mu)\|_Q & = & \sup \{ \| [ \l \Phi(\mu_{e,f}), \phi_{u,v}\r_{\F_Q^*-\F_Q} ] \| : [\phi_{u,v}] \in M_n(\FQA)_{\| \c \|_Q \leq 1} \}\\
 & = & \sup \{ \| [ \l \mu_{e,f}, \phi_{u,v}\r_{\E_Q^*-\E_Q} ] \| : [\phi_{u,v}] \in M_n(\FQA)_{\| \c \|_Q \leq 1} \}\\
 & \geq & \sup \{ \| [ \l \mu_{e,f},   x_{i,j}^* \c x_{k,l}    \r_{\E_Q^*-\E_Q} ] \| : [x_{i,j}^* \c x_{k,l}] \in {\cal S}_n \}\\
 & = &  \sup \{ \| [ \l \mu_{e,f}, \phi_{u,v}\r_{\E_Q^*-\E_Q} ] \| : [\phi_{u,v}] \in M_n(\EQA)_{\| \c \| \leq 1} \}\\
 & = & \| \mu \|_{A^*}
 \eeqs 
 as needed. 
 \end{proof}

 \subsection{The Fourier space $\FqA$}  
 
 Unless stated otherwise, in this section $A$ is a (not necessarily c.c.) Banach algebra and $X$ is a contractive right Banach $A$-module through the action $$q: X \times A \ra X: (x, a) \mapsto x \c a.$$
 By replacing the operator space projective tensor product $\wot$ with the Banach space projective tensor product $\otimes^\gamma$ and dropping the words operator, complete and completely, we obtain the \it Fourier space \rm $\FqA$ and its \it Fourier norm \rm $\| \c \|_q$:  
 $$\FqA = p( X^* \otimes^\gamma X) = p_N(X^* \otimes^\gamma X/N)$$  where $$p:   X^* \otimes^\gamma X \ra A^*: x^* \otimes x \mapsto x^* \c x$$ is the contractive linearization of $q': X^* \times X \ra A^*$, $N = \ker p$ and $\| \c \|_q $ is the Banach space quotient norm inherited from $ X^* \otimes^\gamma X/N $ via the $A$-module isomorphism $p_N: X^* \otimes^\gamma X/N \ra \FqA$, i.e., for $\phi \in \FqA$, 
 \beq   \label{q Norm Formula}     \| \phi \|_q = \inf\{\| \xi \|_\gamma : \xi  \in X^* \otimes^\gamma X \ {\rm and} \ p(\xi) = \phi\}.    \eeq 
  In this situation, we have the following familiar-looking descriptions of $\FqA$ and $\| \c \|_q$. 
 
 \bp \label{FqA Description Prop}    The pair $(\FqA, \| \c \|_q)$ is a  left introverted homogeneous subspace of $A^*$.   If $(x_n^*)$ and $(x_n)$ are sequences in $X^*$ and $X$ respectively, and $\sum \| x_n^*\| \|x_n\| < \infty$, then 
 $\phi = \sum x_n^* \c x_n $ belongs to $\FqA$;  $\sigma(A^*, A)$-convergence of the series implies $\| \c \|_q$-convergence  and $\| \c \|_{A^*}$-convergence. 
Conversely, if $\phi \in \FqA$, then there are sequences  $(x_n^*)$ in $X^*$ and $(x_n)$ in $X$ such that $\sum \| x_n^*\| \|x_n\| < \infty $ and   
 $\phi = \sum x_n^* \c x_n$;  moreover,$$\| \phi \|_q = \inf \left\{ \sum \|x_n^* \| \| x_n\|: \phi = \sum x_n^* \c x_n \right\}. $$
 \ep 
 
 \begin{proof}  The first statement follows from (simpler) versions of the arguments used to establish Theorem \ref{FQA operator introverted Thm}.    Recall that  $\xi \in X^* \otimes^\gamma X$ exactly when there are sequences $(x_n^*)$ in $X^*$ and $(x_n)$ in $X$ such that 
 $$ \sum \|x_n^* \| \| x_n\| < \infty \quad {\rm and } \quad \xi = \sum x_n^* \otimes x_n $$
 with $\| \xi \|_\gamma$ equal to the infimum of $\sum \| x_n^* \| \|x_n\|$ taken over all such representations of $\xi$. Since $p$ is a contractive surjection onto $\FqA$, the validity of the  remaining statements follows from (\ref{q Norm Formula}) via routine arguments. 
 \end{proof}
 
 As in the last subsection, the bilinear map 
 $$ \FqA^* \times X^* \ra X^*: (\mu, x^*) \mapsto \mu \c x^* \ \   where \ \ \l\mu \c x^*, x \r = \l \mu, x^* \c x\r$$ 
 defines a contractive left  $\FqA^*$-module action on $X^*$.  Let  $M_q$ denote the weak$^*$-closure of $\Gamma_q(A)$  in $B(X^*)$, where $\Gamma_q(a) (x^*) = a \c x^*$. Then  $\bg_q: \FqA^* \ra B(X^*)$, where $\bg_q(\mu) (x^*) = \mu \c x^*$, defines a weak$^*$-continuous  isometric algebra isomorphism of $\FqA^*$ onto $M_q$ such that $\bg_q\circ \eta_q = \Gamma_q$. In particular, we have the following proposition. 
 
 \bp \label{FqA* Identification Prop} We can identify the Banach algebra $(\FqA^*, \sq)$ with the weak$^*$-closed operator subalgebra  $M_q$ of $B(X^*)$ via the weak$^*$-homeomorphic  isometric algebra isomorphism 
 $M_q \ra \FqA^*: T \mapsto \mu_T$ defined by $$   \quad \l \mu_T, \sum x_n^* \c x_n \r = \sum \l Tx_n^*, x_n\r \quad whenever \quad \sum \|x_n^*\| \|x_n\| < \infty.$$
 \ep  
  
  Letting $\E_q(A^*) = \E_{\F_q}(A^*)$, the $\|\c \|_{A^*}$-closure of $\FqA$ in $A^*$, the \it Eberlein space \rm associated with $q$, $(\EqA, \| \c \|_{A^*})$, is a left introverted subspace of $A^*$ by Proposition \ref{Introverted Arens product algebras Prop}. 
  
 \br   \label{EqA and FqA Thm Remark}   \rm  The Banach space analogues of the statements found in \ref{FQA and EQA Span Remark}, \ref{EQA Min LDBA Prop}, \ref{HLDBA=FQ*Cor}, \ref{FQA=EQARemark} and \ref{FQA =EQA Prop}  all hold.  \er
  
\br  \label{FQA FqA remark}    \rm  To close this section, assume again that $A$ is a c.c. Banach algebra, $X$ is a c.c. right operator $A$-module via $Q: X \times A \ra X$, and let  $q:X \times A \ra X$ 
denote  the same module action with the operator space structures of $X$ and $A$ ignored. Then we have the containments
$$\FqA \subseteq \FQA \subseteq \EQA = \EqA$$ and, for each $\phi \in \FqA$, $\| \phi\|_Q \leq \| \phi \|_q$. 

To see this, note that because $(\FQA^*, \eta_Q)$ is a HLDBA over $A$ such that $\bg_Q$, when viewed as a mapping of $\FQA^*$ into $B(X^*)$, is a contractive weak$^*$-continuous representation such that $\bg_Q \circ \eta_Q = \Gamma_Q = \Gamma_q$, the Banach space versions of Theorem \ref{EQA Min LDBA Prop}  and Corollary \ref{HLDBA Ordering Cor} imply that  $(\FqA, \| \c \|_q)$ is contained in  $(\FQA, \| \c \|_Q) $ via a contraction.  The equality $\EQA = \EqA$ follows from  Remark \ref{FQA and EQA Span Remark}(a)  and its Banach-space counterpart. \er
  
  \section{Examples} 
  
  In this final section, we will apply our theory to several specific module operations $X \times A \ra X$, thereby recovering -- and often extending -- results concerning (completely) isometric representations  of familiar Banach algebras associated with $A$.  Moreover, by applying Theorem \ref{EQA Min LDBA Prop} to these module actions, we obtain new characterizations of these well-studied Banach algebras.  The identification of some (but not all) of the spaces studied in this section with $\FQA$ and $\EQA$ for some $Q$ could also be achieved by using pre-existing theory in tandem with Corollary \ref{HLDBA=FQ*Cor} and Remark \ref{FQA=EQARemark}. However,  a primary goal in this section is to show how our theory can be applied in a variety of situations to extend and provide new proofs of previously known theorems.  
    
  \subsection{The Fourier spaces $\F_\pi(A^*)$} 
  
  The main objects discussed in this section, $\F_\pi(A^*)$ and its dual,  are new. In the case that $A$ is the group algebra $L^1(G)$ of a locally compact group $G$ and $\prep$ is a continuous unitary representation of $G$, we  recover the Arsac-Fourier spaces $A_\pi$ \cite{Ars}.  
  
  When  $\H$ is a Hilbert space and $\xi \in \H$, we will write $\dot{\xi}$ when viewing $\xi $ as an element of $\ov{\H}$, the conjugate Hilbert space of $\H$.  For $B \in B(\H)$, observe that $\dot{B} \in B(\ov{\H})$ where $\dot{B} (\dot{\xi}) = (B\xi)^{\dot{}}$, and  $B \mapsto \dot{B}$ is a conjugate-linear isometric $*$-isomorphism of $B(\H)$ onto $B(\ov{\H})$. Hence, $\Theta_{\ov{\H}}(B) := (B^*)^{\dot{}} = (\dot{B})^*$ defines a linear isometric anti-homomorphic $*$-isomorphism of $B(\H)$ onto $B(\ov{\H})$. 
  
  Let $A$ be a Banach algebra, $\pi: A \ra B(\H)$ a contractive representation of $A$ on $\H$. Then 
  $\dot{\pi}:= \Theta_{\ov{\H}} \circ \pi: A \ra B(\ov{\H})$ -- i.e., $\dot{\pi}(a) := (\pi(a)^*)^{\dot{}}$ -- is a contractive, linear anti-homomorphic representation of $A$ on $\ov{\H}$.  Therefore, if for $a \in A$ and $\dot{\xi} \in \ov{\H}$ we let 
  $$\dot{\xi} \c a = \dot{\pi}(a)(\dot{\xi}) = (\pi(a)^*(\xi))^{\dot{}}, \quad {\rm then} \quad q_\pi: \ov{\H} \times A \ra \ov{\H}: (\dot{\xi}, a) \mapsto \dot{\xi} \c a$$  defines a contractive right $A$-module action on $\ov{\H}$.  We will use $(\F_\pi(A^*), \| \c \|_\pi)$ to denote the associated Fourier space $(\F_{q_\pi}(A^*), \| \c \|_{q_\pi})$. 
  
  Since $\H$ can be identified with the  dual space of $\oH$ via the linear isometry
  $$\varphi: \H \rightarrow \oH^*: \xi \mapsto \varphi_\xi \ \ \ {\rm given \ by} \ \ \  \l \varphi_\xi, \dot{\eta} \r_{\oH^*-\oH} = \l \dot{\eta}|\dot{\xi} \r_{\oH} = \l \xi| \eta\r_\H,$$ we have 
  $$ q_\pi' : \H \times \oH \ra A^*: (\xi, \dot{\eta}) \mapsto \varphi_\xi \c \dot{\eta}$$ where 
  \beqs  \l \varphi_\xi \c \dot{\eta}, a\r & = &   \l \varphi_\xi, \dot{\eta} \c a \r_{\oH^*-\oH} = \l \dot{\eta} \c a | \dot{\xi} \r_{\oH} \\
  & = & \l (\pi(a)^*(\eta))^{\dot{}} | \dot{\xi} \r_{\oH} = \l \xi | \pi(a)^* \eta \r_\H \\
  & = & \l \pi(a) \xi | \eta\r_\H.
  \eeqs
 Using the coefficient function notation  $\xi*_\pi \eta(a) = \l \pi(a) \xi | \eta\r_\H$ -- cf.  Arsac \cite{Ars} when $A = L^1(G)$, $G$ a locally compact group --  we thus have 
  $$q_\pi': \H \times \oH \ra A^*: (\xi, \dot{\eta}) \mapsto \xpe.$$ 
  
  In the next theorem we will further assume that $A$ is an involutive Banach algebra and $\{ \pi, \H\}$ is a $*$-representation of $A$. (Observe that in this case $\dot{\pi}$ is also a $*$-map.) Of course, a version of this theorem also holds without making these assumptions, the main difference being that $\breve{\pi}$, as defined below, may not be $*$-map and $\F_\pi(A^*)^*$ may fail to be a $W^*$-algebra.

   \bt \label{F_pi(A*) Theorem} Let $A$ be an involutive Banach algebra, $\{\pi, \H\} $ a $*$-representation of $A$.
   \bi  \item[(a)] If $(\xi_n)$, $(\eta_n)$ are sequences in $\H$ and $\sum \| \xi_n\| \|\eta_n\|< \infty$, then $\phi = \sum \xi_n*_\pi \eta_n \in \F_\pi(A^*)$ with $\sigma (A^*, A)$-convergence implying $\| \c \|_\pi$-convergence.  Conversely, if $\phi \in \F_\pi(A^*)$, then there are sequences $(\xi_n)$ and $(\eta_n)$  in $\H$ such that $\sum \| \xi_n\| \|\eta_n\|< \infty$ and $\phi = \sum \xi_n*_\pi \eta_n$; moreover, $$\| \phi \|_\pi =  \inf \left\{ \sum \|\xi_n \| \| \eta_n\|: \phi = \sum \xi_n*_\pi \eta_n \right\}. $$    
 \item[(b)]  $(\F_\pi(A^*), \| \c \|_\pi)$ is a left introverted homogeneous subspace of $A^*$ closed under $\phi \mapsto \phi^*$,  where  $\phi^*(a) = \ov{\phi(a^*)}$. Furthermore, $(\xpe)^* = \epx$ and $\phi \mapsto \phi^*$ is isometric with respect to   $\| \c \|_{A^*} $ and  $ \| \c \|_\pi$.  
 \item[(c)]  With respect to the involution defined by 
 $$\mu^* (\phi)= \overline{\mu(\phi^*)} \qquad \mu \in \F_\pi(A^*)^*, \ \phi \in \F_\pi(A^*),$$
 $(\F_\pi(A^*)^*, \sq)$ is a $W^*$-algebra and $$\breve{\pi}: \F_\pi(A^*)^* \ra B(\H) \quad {\rm given \ by} \quad \l \breve{\pi}(\mu) \xi | \eta\r_\H = \l \mu, \xpe\r_{\F_\pi^*- \F_\pi}$$ 
 is a weak$^*$-homeomorphic isometric $*$-isomorphism of $\F_\pi(A^*)^*$ onto $VN_\pi = \overline{\pi(A)}^{{\rm wk}^*}$, the von Neumann subalgebra of $B(\H)$ generated by $\pi$ on $A$. Moreover, $\breve{\pi} \circ \eta_\pi = \pi$ where $\eta_\pi$ is the canonical homomorphism of $ A$ into $ \F_\pi(A^*)^*$. 
 \item[(d)]  The map $VN_\pi \ra \F_\pi(A^*)^*: T \mapsto \mu_T$ is a weak$^*$-homeomorphic isometric $*$-isomorphism where $$\l \mu_T, \sum \xi_n*_\pi \eta_n\r_{\F_\pi^*-\F_\pi} = \sum \l T\xi_n | \eta_n \r_\H$$ 
 whenever $\sum \xi_n *_\pi \eta_n \in \F_\pi(A^*)$ (with $\sum \|\xi_n\| \|\eta_n\| < \infty$).  In particular, $VN_\pi$ can be identified with the dual of $\F_\pi(A^*)$. 
   \ei
   
   \et 
   \begin{proof}  Part (a) and the statement that $(\F_\pi(A^*), \| \c \|_\pi)$ is a left introverted homogeneous subspace of $A^*$  are immediate consequences of Proposition \ref{FqA Description Prop} and the preceding discussion.  One can quickly verify that $(\xpe)^* = \epx$ and $\| \phi^*\|_{A^*} \leq \| \phi\|_{A^*}$; since $(\phi^*)^* = \phi$, $\|\phi^*\|_{A^*} = \|\phi\|_{A^*}$. If $\phi = \sum \xi_n *_\pi \eta_n \in \F_\pi(A^*)$ with $\sum\|\xi_n\|\|\eta_n\| < \infty$, then part (a) gives $\sum \eta_n *_\pi \xi_n \in \F_\pi(A^*)$ and $\phi^* = (\sum \xi_n*_\pi \eta_n)^* = \sum \eta_n *_\pi \xi_n$ by $\|\c \|_{A^*}$-continuity of $\phi \mapsto \phi^*$. Part (a) now yields $\|\phi^*\|_\pi \leq \| \phi\|_\pi$ which in turn implies that  $\|\phi^*\|_\pi = \| \phi\|_\pi$.  
   
   By Proposition \ref{Introverted Arens product algebras Prop}, $(\F_\pi(A^*)^*, \sq)$ is a Banach algebra and, as noted in the paragraph preceding Proposition \ref{FqA* Identification Prop},   $\breve{\Gamma}_{q_\pi}: \F_\pi(A^*)^* \ra M_{q_\pi}$ is a weak$^*$-homeomorphic isometric algebra isomorphism onto $M_{q_\pi}$, the weak$^*$-closure of $\Gamma_{q_\pi}(A)$ in $B(\oH^*)$. Here
   $$\l \breve{\Gamma}_{q_\pi}(\mu)(\varphi_\xi), \dot{\eta} \r_{\oH^*-\oH} = \l \mu \c \varphi_\xi,  \dot{\eta} \r_{\oH^*-\oH} = \l \mu, \xpe \r$$ and  
   $$\l \Gamma_{q_\pi}(a) (\varphi_\xi), \dot{\eta} \r_{\oH^*-\oH}  = \l a\cdot \varphi_\xi, \dot{\eta} \r_{\oH^*-\oH} = \l \varphi_\xi, \dot{\eta} \c a  \r_{\oH^*-\oH}   = \l \pi(a) \xi | \eta\r_\H.$$ 
   Define an isometric algebra isomorphism, $\kappa$, of $B(\oH^*)$  onto $B(\H)$ by putting $\kappa(B) = \varphi^{-1} \circ B \circ \varphi$. For $B \in B(\oH^*)$ and $\xi, \eta \in \H$, $\varphi_{\kappa(B)\xi} = \varphi \circ \kappa(B) (\xi) = B(\varphi_\xi)$, so 
   $$ \l \kappa(B) \xi | \eta \r_\H = \l \varphi_{\kappa(B)\xi}, \dot{\eta} \r_{\oH^*-\oH} = \l B(\varphi_\xi),  \dot{\eta} \r_{\oH^*-\oH}, $$ from which it follows that $\kappa$ is a weak$^*$-homeomorphism and $\pi = \kappa \circ \Gamma_{q_\pi}$.  Consequently, $$\breve{\pi} := \kappa \circ \breve{\Gamma}_{q_\pi}: \F_\pi(A^*)^* \ra B(\H)$$ is a weak$^*$-continuous isometric isomorphism mapping $\F_\pi(A^*)^*$ onto 
   $$\kappa(M_{q_\pi}) = \kappa\left(\ov{     \Gamma_{q_\pi}(A)}^{{\rm wk}^*}\right) = \ov{   \kappa(\Gamma_{q_\pi}(A))}^{{\rm wk}^*}  = \ov{\pi(A)}^{{\rm wk}^*} = VN_\pi,  $$
  $$\l \breve{\pi}(\mu) \xi | \eta \r_\H = \l \kappa( \breve{\Gamma}_{q_\pi} (\mu)) \xi | \eta\r_\H = \l  \breve{\Gamma}_{q_\pi} (\mu)(\varphi_\xi), \dot{\eta} \r_{\oH^*-\oH} = \l \mu, \xpe\r,$$ 
  and $$\l \breve{\pi}\circ \eta_\pi(a) (\xi) | \eta\r_\H = \l \eta_\pi(a), \xpe \r= \l \pi(a) \xi | \eta \r_\H. $$   
  Let $\mu \in \F_\pi(A^*)^*$. Then $\mu^* \in \F_\pi(A^*)^*$ and 
  $$\l \breve{\pi}(\mu^*) \xi | \eta \r_\H = \ov{ \l \mu, (\xpe)^*\r} = \ov{ \l \mu, \epx \r} = \ov{ \l \breve{\pi}(\mu) \eta | \xi \r_\H   } = \l \breve{\pi}(\mu)^* \xi | \eta \r_\H$$
  so $\breve{\pi}$ is also a $*$-isomorphism onto the von Neumann algebra $VN_\pi$. Hence, $\F_\pi(A^*)^*$ is a $W^*$-algebra, which completes the proof of (c).  Part (d) is a consequence of part (c). 
  \end{proof}

  \br \rm By Theorem \ref{EQA Min LDBA Prop}  and Remark \ref{EqA and FqA Thm Remark}, $(\F_\pi(A^*)^*, \eta_\pi)$ is the minimum HLDBA in $({\cal HLD}(A), \leq)$  [unique HLDBA]  $(\A, \eta_\A)$  over $A$ for which there is a weak$^*$-continuous contractive [isometric] representation $\Gamma_\A : \A \ra B(\H)$ such that $\Gamma_\A \circ \eta_\A = \pi$. 
  \er
  
  Suppose now that $G$ is a locally compact group and $\prep$ is a continuous unitary representation of $G$.  Then $\prep$ determines a non-degenerate $*$-representation of $L^1(G)$ through 
  $$\l \pi(f) \xi | \eta\r_\H = \int \l \pi(s) \xi | \eta\r_\H \,  f(s) \, ds \qquad (f \in L^1(G), \ \xi, \eta \in \H).$$
Using the definition provided above of $\xpe \in L^1(G)^* = L^\infty(G)$, we obtain 
$$ \int \xpe(s) f(s) \, ds = \l \xpe, f\r_{L^\infty - L^1} = \l \pi(f) \xi | \eta \r_\H = \int \l \pi(s) \xi | \eta \r_\H \, f(s) \, ds.$$   
  In $L^\infty(G)$, $\xpe$ therefore equals the continuous coefficient function $s \mapsto \l \pi(s) \xi | \eta\r_\H$  on $G$, as  defined in \cite{Ars}.  
  
  Applying Theorem \ref{F_pi(A*) Theorem}(a) to the $*$-representation $\prep$ of $L^1(G)$, we see from \cite[Theorem 2.2 (ii) and (iii)]{Ars} that our Fourier space $\F_\pi(L^1(G)^*)$ is precisely Arsac's Fourier space $A_\pi$ and our Fourier norm $\| \c \|_\pi$ agrees with the Fourier-Stieltjes algebra norm on $B(G)$ restricted to $A_\pi$.   (In particular, when $\prep$ is the left regular representation $\{\lambda_2, L^2(G)\}$, and  the universal representation   $\{ \omega_G, \H_\omega\}$, of $G$,  $\F_\pi(L^1(G)^*)$ is respectively the Fourier algebra $A(G)$, and the Fourier-Stieltjes algebra $B(G)$.) We also see that the usual identification of $VN_\pi$ with $A_\pi^*$ is a special case of the general result described in Theorem \ref{F_pi(A*) Theorem}(d)  -- which in turn is a special case of Proposition 
  \ref{FqA* Identification Prop}.  Furthermore, Theorem \ref{F_pi(A*) Theorem} identifies $A_\pi$ as a left introverted homogeneous subspace of $L^1(G)^*$, and the product in $VN_\pi$ as an Arens product over $L^1(G)$. 
  
 To reduce the length of this paper, we will postpone the detailed study of the Fourier spaces $\F_\pi(A^*)$ of an involutive Banach algebra $A$  -- and the corresponding Fourier-Stieltjes spaces,  $C^*$-algebras,  von Neumann algebras,  and  Eberlein spaces -- to a subsequent paper.

 \subsection{The Fig\`{a}-Talamanca--Herz spaces $A_p(G)$   } 
  
  Let $G$ be a locally compact group, $1<p< \infty$,  $p'$ the conjugate index satisfying $1/p + 1/p' =1$.  In this subsection, we will observe that $A_p(G)$ is also an example of a Fourier space $\F_q(L^1(G)^*)$, and will recover the identification of $A_p(G)^*$ with the Banach algebra of $p'$-pseudomeasures, $PM_{p'}(G)$.

 Letting $\{ \lambda_p, L^p(G)\}$ be the left regular representation of $G$ on $L^p(G)$ defined by $$\lambda_p(s) \xi(t) = \xi(s^{-1}t) \qquad \xi \in L^p(G), \ s,t \in G,$$  $L^p(G)$ becomes a contractive right Banach $G$-module through $\xi \c s: = \lambda_p(s^{-1}) \xi$. Hence, as noted in Remarks \ref{FQA and EQA Span Remark}(b) and \ref{EqA and FqA Thm Remark}, defining 
 $q: L^p(G) \times L^1(G) \ra L^p(G)$ through the weak integral $\xi \c f = \int \xi \c s \, f(s) \, ds $, $L^p(G)$ is a contractive neo-unital right Banach $L^1(G)$-module such that 
 $$q': L^{p'}(G) \times L^p(G) \ra LUC(G) : (\eta, \xi) \mapsto \eta \c \xi$$
 satisfies 
 $$ \eta \c \xi(s) = \l \eta, \xi \c s \r_{L^{p'}-L^p} = \int \xi  (st) \, \check{\eta}(t^{-1}) \, dt  = \xi * \check{\eta}(s). $$
  Let $(A_p(G), \| \c \|_{A_p} ) = (\F_q(L^1(G)^*), \| \c \|_q)$ and let $PM_p(G)$ denote the weak$^*$-closure in $B(L^p(G))$ of $\lambda_p(L^1(G))$, where $\lambda_p$ is defined on  $L^1(G)$ through the weak integral 
  $ \lambda_p(f) \xi = \int \lambda_p(s) \xi \,  f(s) \, ds $. 
  
  \bc   Let $\phi \in CB(G)$. Then $\phi \in A_p(G)$ if and only if there are sequences $(\xi_n)$ in $L^p(G)$ and $(\eta_n)$ in $L^{p'}(G)$ such that $\sum \|\xi_n\|_p \|\eta_n\|_{p'} <\infty$ and $\phi = \sum \xi_n * \check{\eta}_n$ (with pointwise convergence implying $\|\c\|_{A_p}$ and uniform convergence); moreover 
  $$\| \phi \|_{A_p} = \inf \left\{ \sum \|\xi_n \|_p \| \eta_n\|_{p'}: \phi = \sum \xi_n  * \check{\eta}_n \right\}. $$   Furthermore, 
   $$ \l \breve{\lambda}_{p'}(\mu) (\eta), \xi \r_{L^{p'}-L^p} = \l \mu, \xi * \check{\eta} \r_{A_p^*-A_p} \qquad (\eta \in L^{p'}(G), 
   \ \xi \in L^p(G))$$   
  defines a weak$^*$-homeomorphic isometric algebra isomorphism 
 $\breve{\lambda}_{p'}$ of $ A_p(G)^*$ onto the operator subalgebra  $ PM_{p'}(G)$ of  $B(L^{p'}(G))$.   Thus $PM_{p'}(G)$ and $A_p(G)^*$ can be identified through the pairing 
 $$\l T, \sum \xi_n * \check{\eta}_n\r_{A_p^*-A_p} = \sum \l T \eta_n, \xi_n\r_{L^{p'}-L^p} \ \ 
 {\rm whenever}   \ \ \sum \|\xi_n\|_p \|\eta_n\|_{p'} <\infty.$$  
  \ec 
  
  \begin{proof}  The first statement follows from Proposition \ref{FqA Description Prop} and the above discussion.   As noted in the paragraph preceding Proposition \ref{FqA* Identification Prop}, 
  $$\breve{\lambda}_{p'} (:= \breve{\Gamma}_q) : A_p(G)^* \ra M_q \subseteq B(L^{p'}(G))$$
  is a weak$^*$-homeomorphic isometric algebra isomorphism onto $M_q$, where 
  $$ \l \breve{\lambda}_{p'}(\mu) (\eta), \xi \r_{L^{p'}-L^p} =  \l \mu, \eta \c \xi \r_{A_p^* -A_p} = \l \mu, \xi * \check{\eta} \r_{A_p^*-A_p}$$ and $M_q$ is the weak$^*$-closure of $\breve{\lambda}_{p'} \circ \eta_q (L^1(G)) $ in $B(L^{p'}(G))$; here, as before,  $\eta_q$ is the canonical homomorphism mapping $L^1(G)$ into  $\F_q(L^1(G)^*)^* = A_p(G)^*$. But, for $f \in L^1(G)$, $\xi \in L^p(G)$ and $\eta \in L^{p'}(G)$, 
  \beqs  \l \breve{\lambda}_{p'} ( \eta_q (f)) (\eta), \xi \r_{L^{p'}-L^p}  & = & \l \eta_q (f), \xi * \check{\eta} \r_{A_p^* -A_p}  = \l \xi * \check{\eta}, f\r_{L^\infty-L^1} \\
  & = & \int \int \xi(t) \eta(s^{-1}t) \, dt \, f(s) \, ds = \int \l \lambda_{p'}(s) \eta, \xi \r_{L^{p'}-L^p} f(s) \, ds  \\
  & = & \l \lambda_{p'}(f) \eta, \xi \r_{L^{p'}-L^p}, 
  \eeqs    
  so $\breve{\lambda}_{p'} \circ \eta_q =\lambda_{p'}$. Hence, $M_q =  PM_{p'}(G)$.   
  \end{proof} 
  
  Let $\eta_p$ denote the canonical homomorphism of $L^1(G)$ into $A_p(G)^*$. Observe that by Theorem \ref{EQA Min LDBA Prop}  and Remark \ref{EqA and FqA Thm Remark}, $(A_p(G)^*, \eta_p)$ is the minimum HLDBA in $({\cal HLD}(L^1(G)), \leq)$  [unique HLDBA]  $(\A, \eta_\A)$  over $L^1(G)$ for which there is a weak$^*$-continuous contractive [isometric] representation $\Gamma_\A : \A \ra B(L^{p'}(G))$ such that $\Gamma_\A \circ \eta_\A = \lambda_{p'}$.

  \subsection{The space of left uniformly continuous functionals on $A$, $LUC(A^*)$}
  
 Given a (c.c.) Banach algebra $A$,  recall that $LUC(A^*)$ is the closed linear span of $A^* \c A$,  where $\l a^* \c a,b\r = \l a^*, ab\r$. By recognizing $LUC(A^*)$ as $\EqA = \FqA$ ($\EQA=\FQA$) associated with the right module action of $A$ on itself, with Theorem \ref{LUCRepThm} we establish a (completely) isometric identification of the Banach algebra $LUC(A^*)^*$ with $B_A(A^*)$ ($CB_A(A^*)$), the operator (c.c.) Banach algebra of all (c.b.) $A$-module maps on $A^*$. This extends  each of the corresponding results describing $LUC(A^*)^*$  found  in \cite{Cur-Fig, Lau1, Lau2, Neu-Rua-Spr, Hu-Neu-Rua} for various examples of $A$.  Corollary \ref{LUCasHLBACor} provides new characterizations of $LUC(A^*)^*$.  
  
  Suppose first that $A$ is a (not necessarily c.c.) Banach algebra with a contractive right approximate identity, let $X=A$, and consider the right $A$-module action $$q:X \times A \ra X: (x,a) \mapsto x\c a = xa.$$ Then $$ q': A^* \times A \ra A^*: (a^*,a) \mapsto a^* \c a$$
  where $\l a^* \c a , b\r = \l a^* , ab \r$. Hence, in this case $q'$ is just the usual right $A$-module action on $A^*$.  Since $A$ has a contractive right approximate identity, the Cohen Factorization Theorem  (as stated prior to Corollary \ref{Operator Cohen Factorization Thm})  implies that  $LUC(A^*) = A^* \c A$ and 
  \beq \label{LUC(A*) Unit Ball  Eqn}    LUC(A^*)_{\| \c \| < 1} \subseteq  (A^*)_{\| \c \| \leq 1}  \c A_{\| \c \| \leq 1}  \subseteq LUC(A^*)_{\| \c \| \leq 1}.      \eeq 
 As  $\EqA $ is the $\| \c \|_{A^*}$-closed linear span of $X^* \c X = A^* \c A$, we see that  $\EqA = LUC(A^*)$  in this case. Moreover, (\ref{LUC(A*) Unit Ball  Eqn}) implies that ${\cal S}_1 = \{ x^* \c x: \|x^* \| \leq 1 \ {\rm and} \ \|x\|\leq 1\}$ is $\|\c\|$-dense in the closed unit ball of  $\EqA$, so $$LUC(A^*) = \EqA = \FqA \quad {\rm and} \quad \| \c \|_{A^*} = \| \c \|_q$$
  (see  Remarks \ref{EqA and FqA Thm Remark} and Proposition \ref{FQA =EQA Prop}).

  Now suppose further that $A$ is a c.c. Banach algebra with a contractive right approximate identity. Then $X=A$ is a c.c. right operator $A$-module via $$Q: X \times A \ra X : (x, a)\mapsto x \c a = xa.$$
  Letting $q$ denote the same module action with the operator space structures of $A$ and $X$ ignored, the above discussion and Remarks \ref{FQA FqA remark}  give 
  $$\FqA = \FQA = \EQA = \EqA = LUC(A^*).$$
  Moreover, since $$Q' : A^* \times A \ra A^*: (a^*, a) \mapsto a^* \c a$$
  is a c.c. right operator $A$-module action on $A^*$ and $\EQA = A^* \c A$, by Lemma  \ref{Operator Cohen Factorization Thm},   
  $$M_n(\EQA)_{\| \c \| < 1} \subseteq \left\{[a^*_{k,l} \c a] :  [a^*_{k,l}] \in M_n(A^*)_{\| \c \| \leq 1}, \   a \in  A_{\| \c \|\leq 1} \right\}\subseteq  M_n(\EQA)_{\| \c \|\leq 1}$$ for each positive integer $n$. Since $X=A$, this implies that  ${\cal S}_n$, as defined in Proposition \ref{FQA =EQA Prop}, is  $\|\c \|_{A^*}$-dense in $M_n(\EQA)_{\| \c \| \leq 1}$ for each $n$. Hence, the matrix norms $\| \c \|_Q$ and $\| \c \|_{A^*}$ are equal on $\FQA = \EQA = LUC(A^*)$. 
  
  By Theorem \ref{Completely Isometric Repn Thm}(a), $A^*$ is a c.c. left operator $LUC(A^*)^*$-module via $$LUC(A^*)^* \times A^* \ra A^*: (\mu, a^*) \mapsto \mu \sq a^*$$ where -- using the standard notation --  $\l \mu \sq a^* , a \r = \l \mu, a^* \c a \r$; moreover,  $\eta_{LUC}(a) \sq a^* = a \c a^*$ is the left dual $A$-module action on $A^*$.   Let 
  $$\Gamma: A \ra CB(A^*) \quad {\rm and} \quad \Gamma_{LUC}: LUC(A^*)^* \ra CB(A^*)$$ be defined by 
  $$\Gamma(a)(a^*) = a \c a^* \quad {\rm and} \quad \Gamma_{LUC}(\mu)(a^*) = \mu \sq a^*,$$
  and let $M_{LUC}$ denote the weak$^*$-closure of $\Gamma(A)$ in $CB(A^*)$.   Then by   Theorem \ref{Completely Isometric Repn Thm}(b), $\Gamma_{LUC}$ is a weak$^*$-continuous completely isometric algebra isomorphism of $LUC(A^*)^*$ onto $M_{LUC}$ such that $\Gamma_{LUC} \circ \eta_{LUC} = \Gamma$.   
  
  As an immediate corollary to Theorem \ref{EQA Min LDBA Prop}, we obtain the following new characterizations of $LUC(A^*)^*$. 
  
  \bc \label{LUCasHLBACor} Let $A$ be a c.c. Banach algebra with a contractive right approximate identity. The following statements hold:
  \bi \item[(i)]  $(LUC(A^*)^*, \eta_{LUC})$ is the minimum operator HLDBA [LDBA] $(\A, \eta_\A)$  over $A$  in $({\cal HLD}(A), \leq)$ [$({\cal LD}(A), \leq)$] for which there is a c.c. weak$^*$-continuous representation $\Gamma_\A: \A \ra CB(A^*)$ such that $\Gamma_\A \circ \eta_\A = \Gamma$. 
  
  \item[(ii)]  $(LUC(A^*)^*, \eta_{LUC})$ is the unique operator HLDBA [LDBA] $(\A, \eta_\A)$ over $A$  for which  there is a completely isometric  [faithful] weak$^*$-continuous representation $\Gamma_\A: \A \ra CB(A^*)$ such that $\Gamma_\A \circ \eta_\A = \Gamma$.
   \ei

   \ec
  
  Let $CB_A(A^*)$ denote the c.c Banach algebra comprised of all c.b. right $A$-module maps on $A^*$; so an operator  $T$ in  $ CB(A^*)$   belongs to $CB_A(A^*)$ if $T(a^*\c a ) = T(a^*) \c a $ for each $a^* \in A^*$ and $a \in A$. 
  
  \bt  \label{LUCRepThm} Let $A$ be a c.c. Banach algebra with a contractive right approximate identity.  Then we can identify the c.c. Banach algebra $(LUC(A^*)^*, \sq)$ with the weak$^*$-closed operator subalgebra $CB_A(A^*)$ of $CB(A^*)$ via the weak$^*$-homeomorphic completely isometric algebra isomorphism 
  $$CB_A(A^*) \ra LUC(A^*)^*: T \mapsto \mu_T \ \ {\rm defined \ by } \ \ \l \mu_T, a^* \c a \r = \l Ta^*, a\r.$$
  \et 
  
  \begin{proof}  By Corollary \ref{FQA*=MQ Identification Cor}, $M_{LUC} \ra LUC(A^*)^*: T \mapsto \mu_T$  satisfies the properties described above, so it suffices to show that $M_{LUC}= CB_A(A^*)$. To see this, suppose first that $T \in M_{LUC}$. Then $\mu_T \in LUC(A^*)^*$ and for any $b \in A$, 
  $$\l T(a^* \c a), b\r = \l \mu_T, (a^* \c a) \c b \r = \l \mu_T, a^* \c (a  b) \r  = \l Ta^*, ab\r = \l(Ta^*) \c a, b\r.$$
  Hence, $T \in CB_A(A^*)$.   Assuming that $T\in CB_A(A^*)$, by Proposition \ref{M_Q Description Prop} we can establish that $T \in M_{LUC}$ by showing that $\Lambda_T(a^* \c a) = (Ta^*) \c a$ defines a bounded linear mapping of $LUC(A^*)$ into itself. To see this, first note that if $a^* \c a = b^* \c b$, then 
  $$(Ta^*) \c a = T(a^* \c a ) = T(b^*\c b) = (Tb^*) \c b,$$ so $\Lambda_T$ is a well-defined map. Taking  $\phi_1, \phi_2 \in LUC(A^*) = A^* \c A$, by Theorem 32.23 of \cite{HRII} there are elements $a_1^*, a_2^* \in A^*$ and $a \in A$ such that 
  $\phi_1 = a_1^* \c a$ and $\phi_2 = a_2^*\c a$. Hence, for $\gamma \in \C$, 
  $$\Lambda_T(\gamma \phi_1 + \phi_2) = \Lambda_T((\gamma a_1^* + a_2^*) \c a) = T(\gamma a_1^* +a_2^*) \c a = \gamma (T a_1^*)\c a + (Ta_2^*) \c a = \gamma \Lambda_T(\phi_1) + \Lambda_T(\phi_2),$$  showing that $\Lambda_T$ is linear.  Finally, observe that 
  $\| \Lambda_T(a^* \c a) \| = \| (Ta^*) \c a\| \leq \|T\| \|a^* \| \|a\|$ and, as noted above, 
  $$\S_1 = \{ a^* \c a : \| a^* \| \leq 1 \ {\rm and } \ \| a\| \leq 1\}$$ is $\|\c\|_{A^*}$-dense in $LUC(A^*)$. Hence, $\Lambda_T$ is bounded (with $\|\Lambda_T\| \leq \|T\|$).
  \end{proof}
  
  \br  \rm The same argument shows that if $A$ is a (not necessarily c.c.) Banach algebra with a contractive right approximate identity, then we can identify the Banach algebra $LUC(A^*)^*$ with  $B_A(A^*)$, the weak$^*$-closed operator subalgebra of $B(A^*)$ comprised  of all bounded  right $A$-module maps on $A^*$,  via the weak$^*$-homeomorphic isometric algebra isomorphism defined as in  Theorem \ref{LUCRepThm}. 
  The Banach space version of Corollary \ref{LUCasHLBACor} also holds.  
  \er

    \subsection{Representations of  $LUC(G)^*$ as c.b. maps on $B(\H)$} 
     
   Throughout this section, $G$ is a locally compact group and $\prep$ is a continuous unitary representation of $G$.  Observe that $T(\H)$, the Banach space of  trace-class operators on $\H$ with the trace-class norm $\| \c \|_2$, becomes a right Banach $G$-module through the action $K \c s = \pi(s^{-1}) K \pi(s)$, $(K \in T(\H), \ s \in G)$ (e.g., see \cite[Lemma 2.1]{Bek}) and therefore, as noted in Remark \ref{FQA and EQA Span Remark}(b),  $T(\H)$ is a neo-unital right Banach $L^1(G)$-module via 
   \beqs \label{T(H) L1(G) Module Eqn} Q: T(\H) \times L^1(G) \ra T(\H): (K, f) \mapsto K \c f  = \int K\c s \, f(s) \, ds.   \eeqs     Throughout the remainder of this section, $Q$ will refer to this module action, $s\in G$ and, unless stated otherwise,  $T$ and $K$ are  operators in $B(\H) = T(\H)^*$ and $ T(\H)$, respectively.    
   
   Letting $\xi \otimes \eta^*$ denote the rank-one operator  $(\xi \otimes \eta^*)(\zeta) = \l \zeta| \eta\r \xi$ on $\H$, observe that 
   $(\xi \otimes \eta^*)\c s =(\pi(s^{-1}) \xi) \otimes (\pi(s^{-1})\eta)^*$.  
   
   \bp \label{T(H) Op L1(G) Module Prop}   With respect to  $Q$, $T(\H)$ is a neo-unital c.c. right operator $L^1(G)$-module. 
   \ep 
   
    \begin{proof}  Since $L^1(G)$ has the max operator space structure, $$B(L^1(G), CB(T(\H))) = CB(L^1(G), CB(T(\H))) \cong CB(T(\H) \times L^1(G), T(\H))$$
    \cite[Sections 3.3 and 7.1]{Eff-Rua}.  To establish that $Q$ is a c.c., it  therefore suffices to show that given $f$ in $ L^1(G)$, $Q_f(K)= Q(K,f)$ is c.b. on $T(\H)$ with $\|Q_f \|_{cb} \leq \| f\|_1$. Equivalently, we will show that the dual map $Q_f^*$ is c.b. on $B(\H) = T(\H)^*$ with $\| Q_f^* \|_{cb} \leq \|f\|_1$. To this end, for $T\in M_n(B(\H)) = B(\H^n)$ with $\|T \| \leq 1$ and $\xi, \eta \in \H^n$ with $\| \xi \|, \|\eta\| \leq 1$, we will show that $\ds | \l \left(Q_f^*\right)_n(T) \xi | \eta \r_{\H^n} | \leq \| f\|_1$. 
     We have  
    \beqs | \l \left(Q_f^*\right)_n(T) \xi | \eta \r_{\H^n} | & = & \left|  \sum_{i,j=1}^n  \l Q_f^*(T_{i,j}) \xi_j | \eta_i \r_\H     \right|  =   \left|  \sum_{i,j=1}^n  \l Q_f^*(T_{i,j}),  \xi_j \otimes \eta_i^* \r_{B-T}     \right|  \\
    & = &  \left|  \sum_{i,j=1}^n  \l T_{i,j},   \int (\xi_j \otimes \eta_i^*) \c s \, f(s) \, ds\r _{B-T}     \right | \\
    & = &  \left|  \sum_{i,j=1}^n  \int  \l T_{i,j},    \pi(s^{-1}) \xi_j \otimes (\pi(s^{-1})\eta_i)^* \r _{B-T}  \, f(s) \, ds     \right |  \\
    & = &     \left|   \int \left( \sum_{i,j=1}^n   \l T_{i,j}    \pi(s^{-1}) \xi_j |  \pi(s^{-1})\eta_i \r _{\H} \right)   \, f(s) \, ds     \right |      \\
    & = &     \left|   \int     \l [T_{i,j}]    \left(\pi(s^{-1}) \xi_j\right)_j  |  \left(\pi(s^{-1})\eta_i\right)_i \r _{\H^n}    \, f(s) \, ds     \right | \\     
     & \leq  &     \int    \left| \l T    \left(\pi(s^{-1}) \xi_j\right)_j  |  \left(\pi(s^{-1})\eta_i\right)_i \r _{\H^n} \right|   \,  \left| f(s) \right| \, ds \\
     &\leq & \| f\|_1       
    \eeqs  
    since $\| T\|\leq 1$, $\ds \| \left(\pi(s^{-1}) \xi_j\right)_j  \| = \| \xi \| \leq 1$ and  $\ds \| \left(\pi(s^{-1})\eta_i\right)_i  \| = \| \eta \| \leq 1$.    \end{proof}

    As noted in Remark \ref{FQA and EQA Span Remark}(b), $\F_Q(L^1(G)^*)$ and      $\E_Q(L^1(G)^*)$ are contained in $LUC(G)$ and $T\c K(s) = \l T, K \c s \r_{B-T}$.  Hence, if $K \in T(\H)$ is written as 
    $K = \sum \xi_n \otimes \eta_n^*$, then 
    \beq \label{TK formula}  T \c K(s) = \l T, \sum_{n=1}^\infty  \pi(s^{-1}) \xi_n \otimes (\pi(s^{-1})\eta_n)^* \r_{B-T} =    \sum_{n=1}^\infty \l \pi(s)T \pi(s^{-1}) \xi_n |  \eta_n \r_\H.     \eeq
  Moreover, by Remark  \ref{FQA and EQA Span Remark}(c), 
  $$\mu \mapsto \Gamma(\mu) \in CB(B(\H)) \ \ {\rm where} \ \ \l \Gamma(\mu)(T), K \r = \l \mu, T \c K \r $$ 
  defines a c.c. weak$^*$-continuous  representation of $LUC(G)^*$ as c.b. maps on $B(\H)$. If, as in Section 4.3 of \cite{Sto1}, we employ the notation 
  $$\l\mu, f \r_{LUC^*-LUC} = \int f(s) \, d\mu(s), $$  then we obtain  the following formulation of $\Gamma$ from (\ref{TK formula}): 
  \beqs  \l \Gamma (\mu)(T) \xi | \eta\r_\H = \l \Gamma(\mu)(T) , \xi \otimes \eta^* \r = \l \mu , T \c (\xi \otimes \eta^*) \r = \int \l \pi(s) T \pi(s^{-1}) \xi |  \eta \r_\H  \, d\mu(s).       \eeqs 
  By composing $\Gamma$ on $LUC(G)^*$ with the left-strict--weak$^*$-continuous completely isometric embedding $\Theta$ of $M(G)$ into $LUC(G)^*$ -- note that the restriction of $\Theta$ to $L^1(G)$ is $\eta_{LUC}$ --  we obtain the c.c.  left-strict--weak$^*$-continuous representation 
  $$\Gamma_M: M(G) \ra CB(B(\H)) \ \ {\rm given \ by} \ \  \l \Gamma_M (\mu)(T) \xi | \eta\r_\H = \int \l \pi(s) T \pi(s^{-1}) \xi |  \eta \r_\H  \, d\mu(s) $$ 
  of $M(G)$ as c.b. operators on $B(\H)$; here the left-strict topology on $M(G)$ is taken with respect to the ideal $L^1(G)$ \cite[Lemma 1.2]{Sto}. Observe that  Lemma 5.1 and Theorem 5.2 of  \cite{Fil-Mon-Neu} are contained in these remarks.  
  
  The formulation $$\mu \mapsto \Gamma(\mu) \ \ {\rm where} \ \  \l \Gamma (\mu)(T) \xi | \eta\r_\H = \int \l \pi(s) T \pi(s^{-1}) \xi |  \eta \r_\H  \, d\mu(s) $$ 
  also yields a weak$^*$-continuous completely isometric representation of $\F_Q(L^1(G)^*)^*$ and a c.c.  weak$^*$-continuous faithful representation of       $\E_Q(L^1(G)^*)^*$   as c.b. operators on $B(\H)$ (Theorem \ref{Completely Isometric Repn Thm} and Remarks \ref{FQA and EQA Span Remark}(c)).   When $\prep$ is the left regular representation $\{\lambda_2, L^2(G)\}$ of $G$, we will now show that    $\F_Q(L^1(G)^*)=LUC(G)$ and, as an immediate corollary, will thereby recover the completely isometric  weak$^*$-continuous representation of $LUC(G)^*$ as c.b. mappings on $B(L^2(G))$ due to Neufang \cite{Neu1,Neu2}. (It would be interesting to identify and study $\F_Q(L^1(G)^*)$   and $\E_Q(L^1(G)^*)$    for other continuous unitary representations $\prep$  of $G$.)
  
 For the proof that follows, we note that since  $L^1(G)$ has a (contractive) bounded approximate identity, $LUC(G) = L^\infty(G) \c L^1(G)$; moreover, for $\phi \in L^\infty(G)$ and $f \in L^1(G)$,  $(\phi \c f)(t) = \tilde{f} * \phi(t) = \int \tilde{f}(s) \phi(s^{-1}t) \, ds $ where $\tilde{f}(s)  = \Delta(s^{-1}) f(s^{-1})$.

  \bt \label{LUC(G)=FQ Thm}  Let $Q: T(\H) \times L^1(G) \ra T(\H)$ be taken with respect to $\{ \lambda_2, L^2(G)\}$, $q$ the same module action with the operator space structures ignored. Then 
  \beq  \label{LUC(G)=FQ Thm Eq}  LUC(G) = \E_Q(L^1(G)^*) = \F_Q(L^1(G)^*) = \F_q(L^1(G)^*),   \eeq
  the matrix norms $\| \c \|_Q$ and $\|\c \|_{L^1(G)^*}$ agree on this common space, and $\| \c \|_q=\| \c \|_{L^1(G)^*}$ on $LUC(G)$. 
  \et   
    
  \begin{proof}   We have already observed that each of the spaces displayed in (\ref{LUC(G)=FQ Thm Eq}) is contained in $LUC(G)$. Let $\psi = [\psi_{i,j}] \in M_n(\E_Q(L^1(G)^*))$ with $\| \psi\|_{L^1(G)^*}  < 1$.  To see that $ \E_Q(L^1(G)^*) = \F_Q(L^1(G)^*) $ and the matrix norms  $\| \c \|_Q$ and $\|\c \|_{L^1(G)^*}$ agree, by Proposition  \ref{FQA =EQA Prop} it suffices to find $T= [T_{i,j}] \in    M_n(B(\H))_{\| \c \| \leq 1}$ and $K \in T(\H)_{\| \c \|\leq 1}$ such that $[\psi_{i,j}] = [T_{i,j} \c K]$. 
  
  Since $(\phi, f) \mapsto \phi \c f = \tilde{f} * \phi$ defines a c.c. right operator (dual) $L^1(G)$-module action on $L^\infty(G)$, $LUC(G) = L^\infty(G) \c L^1(G)$ and $\psi \in M_n(LUC(G))_{\| \c \| <1}$, it follows from the Cohen factorization theorem -- specifically Lemma \ref{Operator Cohen Factorization Thm} -- that $[\psi_{i,j}] = [\phi_{i,j} \c f]$  for some $[\phi_{i,j}]$ in $M_n(L^\infty(G))_{\| \c \| \leq 1}$  and $f \in L^1(G)_{\| \c \| \leq 1}$.  The map $L^\infty(G) \ra B(\H): \phi \mapsto M_{\check{\phi}}$ where $\check{\phi}(s) = \phi(s^{-1}) $ and $M_{\check{\phi}}$ is the multiplication operator on $\H=L^2(G)$ by $\check{\phi}$ is a composition of $*$-isomorphisms of $C^*$-algebras, and is therefore a complete isometry \cite[p. 26]{Eff-Rua}.  Letting $T_{i,j} = M_{\check{\phi}_{i,j}}$, we therefore have $T = [T_{i,j}] \in M_n(B(\H)) $ with $\| T \|\leq 1$. Following the proof of \cite[Theorem 5.3]{Fil-Mon-Neu}, we now define $\xi, \eta \in \H$ by letting $\xi = |\tilde{f}|^{1/2}$, $\eta(t) = \tilde{f}(t) / |  \tilde{f}|^{1/2}(t)$ if $       \tilde{f} (t) \neq 0$ and $\eta(t) = 0$ otherwise. Then $K = \xi \otimes \eta^* \in T(\H)$ with $\|K \|_1 = \| \xi \|_2 \|\eta\|_2 = \| \tilde{f} \|_1 = \| f\|_1 \leq 1$ and for $t\in G$, 
  \beqs 
  T_{i,j} \c K(t) &  = & \l T_{i,j} \lambda_2(t^{-1}) \xi | \lambda_2(t^{-1})\eta \r = \int \phi_{i,j}(s^{-1}) \xi(ts) \ov{\eta(ts)} \, ds \\
  & = & \int \phi_{i,j}(s^{-1}) \tilde{f} (ts)  \, ds  = \int \tilde{f}(s) \phi_{i,j}(s^{-1}t) \, ds \\
 & = & \tilde{f} * \phi_{i,j}(t) = (\phi_{i,j} \c f) (t) = \psi_{i,j}(t).  \eeqs 
  Thus, $[\psi_{i,j}] = [T_{i,j} \c K]$ with $\|T\| \leq 1$ and $\| K \|_1 \leq 1$, as needed.       
   
   Observe that the above argument shows that any $\psi \in LUC(G)$ with $\| \psi \| <1$ can be written as $\psi = T \c K \in \F_Q(L^1(G)^*)$ for some $T \in B(\H)_{\| \c \|  \leq 1}$ and $K \in T(\H)_{\| \c \| \leq 1}$. Hence,  $LUC(G) = \F_Q(L^1(G)^*) = \E_Q(L^1(G)^*)= \E_q(L^1(G)^*)$. Moreover, $\E_q(L^1(G)^*)= \F_q(L^1(G)^*)$ and  $\| \c \|_q = \| \c \|_{L^1(G)^*}$ by the Banach space version of Proposition \ref{FQA =EQA Prop} (see Remarks \ref{EqA and FqA Thm Remark}).    
   \end{proof}  
    
   The following is an immediate  corollary to Theorems \ref{LUC(G)=FQ Thm} and \ref{EQA Min LDBA Prop}.   Part (i) was first proved by Neufang \cite{Neu1, Neu2}; parts (ii) and (iii) are new. 
  
  \bc  \label{Neufang LUC(G)* Representation Cor}   The following statements hold: 
  \bi \item[(i)] The map  
    $$\Gamma: LUC(G)^* \ra CB(B(L^2(G))) \ \ {\rm given \ by} \ \  \l \Gamma (\mu)(T) \xi | \eta\r = \int \l \lambda_2(s) T \lambda_2(s^{-1}) \xi |  \eta \r  \, d\mu(s) $$  
    defines a completely isometric  weak$^*$-continuous representation of $LUC(G)^*$ as c.b. mappings on $B(L^2(G))$.
    \item[(ii)]  $(LUC(G)^*, \eta_{LUC})$ is the unique operator HLDBA [LDBA] $(\A, \eta_\A)$ over $A$  for which  there is a completely isometric  [c.c. faithful] weak$^*$-continuous representation $\Gamma_\A: \A \ra CB(B(L^2(G)))$ such that $\Gamma_\A \circ \eta_\A(f)  = \Gamma(f)$ for each $f \in L^1(G)$;     
\item[(iii)]  $(LUC(G)^*, \eta_{LUC})$ is the minimum operator HLDBA [LDBA] $(\A, \eta_\A)$  over $L^1(G)$  in $({\cal HLD}(L^1(G)), \leq)$ [$({\cal LD}(L^1(G)), \leq)$] for which there is a c.c. weak$^*$-continuous representation $\Gamma_\A: \A \ra CB(B(L^2(G)))$ such that $\Gamma_\A \circ \eta_\A(f)  = \Gamma(f)$ for each $f \in L^1(G)$. 
    \ei    \ec
    
 Composing $\Gamma$ --  as defined in Corollary \ref{Neufang LUC(G)* Representation Cor} -- with the embedding  $\Theta:M(G) \hookrightarrow LUC(G)^*$, yields the completely isometric representation of $M(G)$ as c.b. operators on $B(L^2(G))$ due to St$\o$rmer  in the abelian case \cite{Stor},  Ghahramani in the general isometric form \cite{Gha},  and Neufang, Ruan and Spronk  in the completely isometric form \cite{Neu-Rua-Spr}.    
 
 \subsection{Other examples} 
 
 There are many other examples of operator HLDBAs  $(\A, \eta_\A)$ for which there is a completely isometric weak$^*$-continuous representation $\Gamma_\A: \A \ra CB(X^*)$ such that $\Gamma_\A \circ \eta_\A$ maps $A$ into $CB^\sigma(X^*)$. 
 
 For instance, if $A$ is a (c.c.) Banach algebra and $\SA$ is a closed introverted subspace of $A^*$ contained in $WAP(A^*)$ -- the space of weakly almost periodic functionals on $A^*$ -- then $\SA^*$ is a (operator) dual Banach algebra over $A$ and therefore, by \cite[Corollary 3.8]{Daws} (the main result of \cite{Uyg}),  there is a reflexive (operator) space $E$ and a (completely) isometric weak$^*$-continuous representation $\Gamma_\S$ of $\SA^*$ into  $B(E)= B^\sigma(E)$ ($CB(E)= CB^\sigma(E)$).  When $A = L^1(\mathbb{G})$ for some locally compact quantum group $\mathbb{G}$, the authors of \cite{Hu-Neu-Rua} studied several examples of closed left introverted subspaces of $A^* = L^\infty(\mathbb{G})$ for which there is an operator space $X$ and   a completely isometric weak$^*$-continuous representation $\Gamma_\S:\SA^* \ra CB(X^*)$ such that $\Gamma_\S \circ \eta_\S$ maps  $L^1(\mathbb{G})$ into $CB^\sigma(X^*)$. 
 
 By Corollary \ref{HLDBA=FQ*Cor} and Remarks \ref{FQA=EQARemark}, for each of these left introverted spaces there is a c.c. right operator $A$-module action $Q$ on $X$ such that $\SA= \EQA = \FQA$, the matrix norms $\| \c \|_{A^*}$ and $\| \c \|_Q$ are equal on this common space and, therefore, the operator HLDBAs $(\SA^*, \eta_\S)$, $(\FQA^*, \eta_Q)$ and $(\EQA^*, \eta_\E)$ coincide. (Theorem \ref{EQA Min LDBA Prop} thus provides new characterizations of these LDBAs over $A$.)  It would be interesting to  try establishing some of these statements  directly as we did with the other examples in this section, and thereby obtain new proofs of these representation theorems. 
  \bigskip 
  
  \noindent {\bf Acknowledgements:} The author is grateful to the anonymous referee whose comments have improved the exposition of this paper.

\noindent {\sc Department of Mathematics and Statistics, University
of Winnipeg, 515 Portage Avenue, Winnipeg, MB, Canada, R3B 2E9 }

\noindent email: {\tt r.stokke@uwinnipeg.ca}


\begin{thebibliography}{99}

\bibitem{Are}  R. Arens,  The adjoint of a bilinear operation, {\it Proc. Amer. Math. Soc.} 2, (1951), 839-848. 

\bibitem{Ars} G. Arsac,  Sur l'espace de Banach engendr$\acute{\rm{e}}$
par les coefficients d'une repr$\acute{\rm{e}}$sentation unitaire,
{\it Publ. D$\acute{e}$p. Math. (Lyon)} 13 (1976), 1-101.



\bibitem {Bek} M.E.B.  Bekka,   Amenable unitary representations of locally compact groups, {\it Invent. Math. }100, 1990, 383-401.

\bibitem{Ber-Jun-Mil} J.F. Berglund, H. Junghenn and P.  Milnes,  {\it Analysis on semigroups: Function spaces, compactifications, representations,}  Canadian Mathematical Society Series of Monographs and Advanced Texts,  John Wiley \& Sons, Inc., New York, 1989.

\bibitem{Con} J. Conway, \it A course in functional analysis,  second edition, \rm Graduate Texts in Mathematics, 96, Springer-Verlag, New York, 1990. 

\bibitem{Cur-Fig} P. C. Curtis Jr. and A. Figa-Talamanca, Factorization theorems for Banach algebras, {\it Function Algebras (Proc. Internat. Sympos. on Function Algebras, Tulane Univ., 1965)}, Scott-Foresman, Chicago, I11., 1966, pp. 169-185. 

\bibitem{Dal} H.G. Dales, {\it Banach algebras and automatic
continuity}, London Math. Soc. Monographs, Volume 24, Clarendon
Press, Oxford, (2000).


\bibitem{Dal-Lau}  H.G. Dales and A. T.-M. Lau,  {\it The second duals of Beurling algebras},  Mem. Amer. Math. Soc. 177 (2005), no. 836, vi+191 pp.
 
 \bibitem{Daws}  M. Daws, Dual Banach algebras: representations and injectivity, {\it Studia Math.} 178 (2007), no. 3, 231-275.  

\bibitem{Eff-Rua} E.G. Effros and Z.-J. Ruan, {\it Operator Spaces}, Oxford
University Press, 2000.

 \bibitem{Fil-Mon-Neu} M. Filali, M. Neufang and M. Sangani Monfared, Representations of Banach algebras subordinate to topologically introverted spaces, {\it Trans. Amer. Math. Soc. } 367 (2015), no. 11, 8033-8050.
 
 \bibitem{Gha} F. Ghahramani, Isometric representation of $M(G)$ on $B(H)$,  \it Glasgow Math. J. \rm  23 (1982), no. 2, 119-122.  
 
 \bibitem{Hu-Neu-Rua}  Z. Hu, M. Neufang and Z.-J. Ruan, Multipliers on a new class of Banach algebras, locally compact quantum groups, and topological centres, {\it Proc. Lond. Math. Soc.}  (3) 100 (2010), no. 2, 429-458. 
 
\bibitem{HRII} E. Hewitt and K.A. Ross, {\it Abstract Harmonic Analysis II}, Grundelehernder mathematischen Wissenschasften 152, Springer, Berlin, 1970.

\bibitem{Lau1}  A.T.-M. Lau,  Uniformly continuous functionals on the Fourier algebra of any locally compact group, {\it Trans. Amer. Math.} Soc. 251 (1979), 39-59. 

\bibitem{Lau2}  A.T.-M. Lau,  Uniformly continuous functionals on Banach algebras, {\it  Colloq. Math.} 51 (1987), 195-205. 

\bibitem{Neu1} M. Neufang, {\it Abstrakte Harmonische Analyse und Modulhomomorphismen $\ddot{u}$ber von Neumann-Algebren}, Ph.D. Thesis, UniversitŠt des Saarlandes (2000).

\bibitem{Neu2} M. Neufang, Isometric representation of convolution algebras as completely bounded module homomorphisms and a characterization of the measure algebra, unpublished manuscript.

\bibitem{Neu-Rua-Spr}   M. Neufang, Z.-J. Ruan, and N.  Spronk,  Completely isometric representations of $M_{cb}A(G)$ and  $UCB(\widehat{G} )$, {\it Trans. Amer. Math. Soc.}  360 (2008), no. 3, 1133-1161. 


\bibitem{Pis} G. Pisier, {\it Introduction to operator space theory},  London Mathematical Society Lecture Note Series, 294,  Cambridge University Press, Cambridge, 2003.

\bibitem{Rua} Z.-J. Ruan,  The operator amenability of A(G), {\it Amer. J. Math.} 117 (1995), no. 6, 1449-1474.

\bibitem{Spr-Sto} N. Spronk and R. Stokke,  Matrix coefficients of unitary representations and associated compactifications of locally compact groups, {\it Indiana Univ. Math. J. }   62 (2013), no. 1, 99-148.

 \bibitem{Sto} R. Stokke, Homomorphisms of convolution algebras, {\it J. Funct. Anal.} 261 (2011), no. 12, 3665-3695.
 
 \bibitem{Sto1} R. Stokke, Amenability and modules for Arens product algebras {\it  Q. J. Math.}  66 (2015), no. 1, 295-321.
  
  
  

\bibitem{Stor} E. St${\o}$rmer,  Regular abelian Banach algebras of linear maps of operator algebras, \it J. Funct. Anal. \rm  37 (1980), no. 3, 331-373.

\bibitem{Uyg} F. Uygul, A representation theorem for completely contractive dual Banach algebras, {\it J. Operator Theory} 62 (2009), no. 2, 327-340. 


 \end{thebibliography}
\end{document}